\documentclass{amsart}
\usepackage{amsmath,amsthm,amssymb,amsfonts,ifpdf,bbm,verbatim}
\usepackage{mathabx,esint}
\usepackage{enumitem}
\usepackage{graphicx}
\usepackage[usenames,dvipsnames]{color}
\usepackage{tensor,calc}

\ifpdf
\usepackage[pdftex,pdfstartview=FitH,pdfborderstyle={/S/B/W 1},%
colorlinks=true, linkcolor=blue, urlcolor=blue, citecolor=blue,%
pagebackref=true]{hyperref}
\else
\usepackage[dvips]{hyperref}
\fi
\numberwithin{equation}{section}

\usepackage[shortalphabetic,abbrev,msc-links,nobysame]{amsrefs}

\newtheorem {Theorem}    {Theorem}[section]
\newenvironment{Theorem*}
  {\Theorem}
  {\endTheorem}

\newtheorem {Corollary}  [Theorem]    {Corollary}
\newtheorem {Proposition}[Theorem]    {Proposition}

\newtheorem {Conjecture}{Conjecture}

\theoremstyle{definition}
\newtheorem {Definition} [Theorem]    {Definition}

\newtheorem {Exercise}   {Exercise}
\newenvironment{Exercise*}
  {\Exercise}
  {\endExercise}
\newtheorem {Question}   [Theorem]    {Question}
\newcounter{AbcT}

\numberwithin{equation}{section}

\newcommand {\Heads}[1]   {\smallskip\pagebreak[1]\noindent{\bf #1{\hskip 0.2cm}}}
\newcommand {\Head}[1]    {\Heads{#1:}}

\newcommand {\C} {{\mathbb C}}

\newcommand {\N} {{\mathbb N}}
\newcommand {\Q} {{\mathbb Q}}
\newcommand {\R} {{\mathbb R}}

\newcommand {\T} {{\mathbb T}}
\renewcommand {\H} {{\mathbb H}}
\newcommand {\Z} {{\mathbb Z}}

\renewcommand{\liminf}{\varliminf}
\renewcommand{\limsup}{\varlimsup}
\DeclareMathOperator{\supp}{supp}
\DeclareMathOperator{\stab}{stab}
\DeclareMathOperator{\cl}{cl}

\DeclareMathOperator{\SL}{SL}
\DeclareMathOperator{\PSL}{PSL}

\DeclareMathOperator{\GL}{GL}

\DeclareMathOperator{\PGL}{PGL}
\DeclareMathOperator{\SO}{SO}

\DeclareMathOperator{\SU}{SU}

\newcommand {\IGNORE}[1]  {}

\renewcommand {\setminus}       {\smallsetminus}

\newcommand {\absolute}[1] {\left| {#1} \right|}
\newcommand {\norm}[1] {\left\| {#1} \right\|}

\newcommand{\GG}{\mathbb{G}}

\newcommand\vol{\operatorname{vol}}

\newcommand\cA{\mathcal{A}}

\DeclareMathOperator{\trace}{trace}
\DeclareMathOperator{\grass}{gr}
\DeclareMathOperator{\size}{size}
\DeclareMathOperator{\disc}{disc}

\newtheorem*{ConjectureA}{Conjecture}

\setlist[enumerate,1]{leftmargin=*}
\setlist[itemize,1]{leftmargin=*}
\makeatletter
\@namedef{subjclassname@2020}{%
  \textup{2020} Mathematics Subject Classification}
\makeatother

\begin{document}

	\title[Rigity properties of higher rank diagonalizable actions]{Recent progress on rigity properties of higher rank diagonalizable actions and applications}
	\author[E. Lindenstrauss ]{Elon Lindenstrauss}
	\address{The Einstein Institute of Mathematics\\
	Edmond J. Safra Campus, Givat Ram, The Hebrew University of Jerusalem
	Jerusalem, 91904, Israel}
	\email{elon@math.huji.ac.il}
	\thanks{E.~L.~acknowledges the support of ERC-2018-ADG project HomDyn and ISF grant 891/15.}
	\date{\today}
	\begin{abstract} The rigidity propeties of higher rank diagonalizable actions is a major theme in homogenous dynamics, with origins in work of Cassels and Swinnerton-Dyer in the 1950s and Furstenberg. We survey both results and conjectures regarding such actions, with emphasize on the applications of these results towards understanding the distribution of integer points on varieties, quantum unique ergodicity, and Diophantine approximations. 
	\end{abstract}
	\dedicatory{Dedicated to G.A.~Margulis}
	\maketitle

\section{Introduction}\label{introduction}
The extensive theory of actions of unipotent groups on homogeneous spaces, to which G.\ A.\ Margulis made many pioneering contributions, gives very satisfactory qualitative (if not yet quantitative) understanding of these actions, with numerous and profound applications. The current state-of-the-art regarding actions of diagonalizable groups is much less satisfactory, and indeed for most natural questions we only have partial results regarding the dynamics. Fortunately, these partial results already have fairly wide applicability. It is the purpose of the survey to present some of the rigidity results regarding such actions as well as their applications.

The motivation to studying rigidity properties of higher rank diagonal actions comes from two different directions. One of these is from the geometry of numbers: the program, initiated by Minkowski, of using lattices in Euclidean spaces and their generalizations for understanding number theoretic questions. We shall make in the survey a distinction between \emph{arithmetic questions}, that is to say properties of integer points (or more generally rational or algebraic point), such as counting and distribution properties of integer points on varieties, and questions in \emph{Diophantine approximations}, such as how well can a point to a given variety be approximated by integer or rational points. Both kinds of applications were already prominently present in the geometry of numbers since its inception by Minkowski.

A classical problem in the geometry of numbers is the study of the set of values attained at integer points by a homogeneous form $F$ of degree $d$ obtained by taking the products of $d$-linear forms in $d$-variables, i.e. one considers forms
\begin{equation*}
F (x _ 1, \dots, x _ d) = \prod_ {i = 1 } ^ d l_i(x_1, \dots, x _ d)
\end{equation*}
where $l_1, \dots,l_d$ are $d$ linearly independent linear forms\footnote{When considering a product of $d$ linear forms in $d$ variables the forms will be implicitly assumed to be linearly independent even if this is not explicitly stated.}, and investigate the values attained by $F$ for $\mathbf x = (x _ 1, \dots, x _ d) \in \Z ^ d$. For instance, one may study the quantity
\[ \nu _ F = \inf_ {\mathbf x \in \Z ^ d \setminus \left\{ 0 \right\} } \absolute {F(\mathbf x)}
.\]
If we present the coefficients of the linear forms $l _ i$ in a $d \times d$-matrix $g$ (one row for each linear form) then the map $F(g)$ assigning a product of $d$ linear forms $F$ to a $d\times d $-matrix $g$ is left invariant under the action of the $d -1$-dimensional diagonal subgroup $A<\SL (d, \R)$ whereas the map $F \mapsto \nu _ F$ is invariant under composition of $F$ by an element of $\GL (d, \Z)$ (in the geometry of numbers literature, two forms which are the same up to the action of $\GL (d, \Z)$ are said to be \emph{equivalent}\footnote{Sometimes one makes a distinction between forms that are the same up to composition by an element of  $\SL (d, \Z)$, which are said to be properly equivalent, and forms that are the same under the action of the slightly bigger group $\GL (d, \Z)$ which are only said to be improperly equivalent.}).
Thus we may view $g \mapsto \nu_{F(g)}$ as either a
(left) $A$-invariant function on $\GL (n, \R) / \GL (n, \Z)$ or a (right) $\GL (n, \Z)$-invariant function on $A \backslash \GL (n, \R)$.
It is convenient to use the normalized quantity  $\bar {\nu} _ {F (g)} = \nu _ {F(g)}/\absolute {\det g}$ which is a well defined function on $\PGL (d, \R)$, left invariant under $A$ and right invariant under $\PGL (d, \Z)$.

Already the case $d = 2$ is of some interest and was quite extensively studied \cite[\S II]{Cassels-Diophantine-book}. In this case the possible forms $F$ considered  are simply the set of  non-degenerate indefinite quadratic forms in two variables. 
For any product of $2$ linear forms in $2$ variables $F$, the value of $\bar {\nu}_F$ is $\leq \sqrt 5$, with equality if and only if $F$ is equivalent (up to a multiplicative scalar, and the action of $\GL (2,\Z)$) to $F(x,y)=x^2-xy-y^2$. Up to the same degrees of freedom, Markoff constructed a complete list of (countably many) binary form with $ \bar {\nu}_F > \frac 13$ but there are uncountably many such indefinite binary forms with $ \bar {\nu}_F = \frac 13$ and a set of full Hausdorff dimension of forms with $\bar {\nu}_F>0$.
Cassels and Swinnerton-Dyer investigated the possible values of $\bar \nu _ F$ for forms that are a product of three linear forms in three variables \cite{Cassels-Swinnerton-Dyer}. They discovered that integral forms of this type satisfy a very strong isolation result, much stronger than the analogous isolation result of Remak and Rogers for product of two linear forms in two variables. We emphasize that an integral form that is a product of $d$ linear forms in $d$ variables need not be presentable as a product of $d$ integral linear forms in $d$ variables. This led them to make the following remarkable conjecture (to be precise, Cassels and Swinnerton Dyer state this conjecture in their paper for $d=3$, but it is clear that they realized a similar phenomenon should hold for higher $d$; cf.\ also the much later remark in Swinnerton-Dyer's book \cite[p.20]{Swinnerton-Dyer-book}):
\begin{Conjecture} [Cassels and Swinnerton-Dyer\cite{Cassels-Swinnerton-Dyer}]\label{C-SD conjecture} Let $d \geq 3$. Any form $F$ which is a product of $d$ linear forms in $d$ variables which is not proportional to a form with integral coefficients has $\bar \nu _ F=0$.
\end{Conjecture}
This farsighted paper of Cassels and Swinnerton-Dyer, and in particular the above conjecture, was highlighted by Margulis in \cite{Margulis-Oppenheim-conjecture}. Stated in terms of the homogeneous space $\PGL (d, \R) / \PGL (d, \Z)$ this conjecture is equivalent to the following:
\begin{ConjectureA}[\cite{Cassels-Swinnerton-Dyer,Margulis-Oppenheim-conjecture}]
Let $d \geq 3$. Any orbit of the diagonal group $A$ in $\PGL (d, \R) / \PGL (d, \Z)$ is either unbounded or periodic.
\end{ConjectureA}

\noindent
Here and throughout we say that an orbit $L.x$ of a locally compact group $L$ on a space $X$ is periodic if the stabilizer of $x$ is a lattice in $L$, i.e. discrete and of finite covolume. Cassels and Swinnerton-Dyer show in \cite{Cassels-Swinnerton-Dyer} that Conjecture~\ref{C-SD conjecture} implies a conjecture of Littlewood from c. 1930:
\begin{Conjecture} [Littlewood]\label{Littlewood Conjecture}
For any $\alpha, \beta \in \R$, it holds that \[\inf_ {n>0} n \norm {n \alpha} \norm {n \beta} = 0.\]
\end{Conjecture}

\noindent 
Here we use for $x \in \R$ the somewhat unfortunate but customary notation $\norm x = \min_ {n \in \Z} \absolute {x - n}$.

A second historical motivation comes from ergodic theory, namely the work of Furstenberg on ``transversality'' of the $\times a$ and $\times b$ maps on $\T=\R / \Z$ for $a$ and $b$ multiplicatively independent. Recall that two integers $a$ and $b$ are said to be multiplicatively independent if they are not both powers of this same integer, i.e. if $\log a / \log b \not\in \Q$.
In his landmark paper \cite{Furstenberg-disjointness-1967}, Furstenberg proved the following theorem:

\begin{Theorem} [Furstenberg \cite{Furstenberg-disjointness-1967}] \label{Furstenberg's theorem}
Let  $X$ be a closed subset of $\T$ invariant under the action of the multiplicative semigroup $S_{a,b,}$ generated by two multiplicatively independent integers $a$ and $b$, that is to say $s.x \in X$ for any $s \in S$ and $x \in X$. Then $X$ is either finite or $X=\T$.
\end{Theorem}
In this paper Furstenberg introduces the notion of joinings and the related notion of disjointness of dynamical systems, a notion that would be important for us later on in this survey, and deduced Theorem~\ref{Furstenberg's theorem} from a particular disjointness principle, one of several enunciated in the paper.

He also presented the following highly influential conjecture that is still open, a natural analogue to Theorem~\ref{Furstenberg's theorem} in the measure preserving category:
\begin{Conjecture} [Furstenberg, c.\ 1967]\label{Furstenberg's conjecture} Let $S_{a,b} \subset \N$ be a semigroup generated by two multiplicatively independent integers as above, and let $\mu$ be a $S_{a,b}$-invariant and ergodic probability measure on $\T$. Then either $\mu$ is finitely supported or $\mu$ is the Lebesgue measure $m_{\T}$ on $\T$.
\end{Conjecture}

\noindent
In this survey will use Greek letters to denote unknown probability measures, and $m$ (often decorated with subscripts) to denote a ``canonical'' probability measures such as Lebesgue measure or Haar measure. We stress that $\mu$ being $S _ {a, b}$-ergodic does not imply it is ergodic under the action generated by multiplication by $a$ or by $b$ --- only that any measurable subset $X \subset \T $ which is $S _ {a, b}$-invariant has either $\mu (X) = 0$ or $\mu (\T \setminus X) = 0$. 

Dealing with semigroup actions is somewhat awkward; this is easily remedied though: it is easy to see that Conjecture \ref{Furstenberg's conjecture} is equivalent to the classifying the $\left\{ a^k b^l:k,l \in \Z\right \}$-invariant and ergodic probability measures on $\prod_ {p | ab \ \text{prime}} \Q _ p /\Lambda$ with $\Lambda = \Z[1/ab]$ diagonally embedded in $\prod_ {p | ab} \Q _ p$.

\label{Rudolph theorem mention}
An important insight of Rudolph \cite{Rudolph-2-and-3}, building upon prior work of Lyons \cite{Lyons-2-and-3}, is that entropy plays an important role in understanding this measure classification question. Specifically, Rudolph proved for $a,b$ relatively prime that Lebesgue measure is the only $S _ {a,b}$-invariant and ergodic probability measure on $\T$ so that its entropy with respect to at least one  element of $S _ {a,b}$ is positive. This was extended to the more general multiplicative independent case by Johnson \cite{Johnson-invariant-measures}.
We now have quite a few other proofs of Rudolph's theorem (e.g. \cite{Feldman-generalization,Host-normal-numbers,Hochman-times-m} to name a few) that seem to me quite different, though all rely very heavily on the positive entropy assumption.

The key feature of the rigidity of higher rank abelian groups such as the action of $S _ {a, b}$ on $\T$ is that the rigidity does not come from the action of an individual element. For any $s \in S _ {a, b}$ there are uncountably many $s$-invariant and ergodic probability measures on $\T$ with any entropy in $[0,\log s]$ as well as uncountably many $s$-invariant closed subset of $\T$ of any Hausdorff dimension in the range $[0,1]$, though we do mention one important restriction: Lebesgue measure $m _ \T$ is the unique $s$-invariant measure on $\T$ with entropy $\log s$, and $\T$ is only $s$-invariant closed subset of $\T$ of Hausdorff dimension 1.

Furstenberg presented the $\times a$, $\times b$-problem as a special instance of a more general problem, and indeed the type of phenomena pointed out by Furstenberg exists also in the action of $A$ on $\PGL (d, \R) / \PGL (d, \Z)$ and in many other high dimensional diagonal actions.

The key feature of rigidity of higher rank diagonalizable actions --- rigidity of the action as a whole while no rigidity for the action of individual elements --- is in contrast to the rigidity properties of unipotent groups and more generally actions of groups generated by unipotents. 

\begin{Definition} Let $G$ be a locally compact group, $\Gamma < G$ a closed subgroup. A measure $\mu$ on $G / \Gamma$ is said to be \emph{homogeneous} if it is supported on a single orbit of its stabilizer $\stab _{G} \mu = \left\{ g \in G: g. \mu = \mu\right \}$.
\end{Definition}

A landmark result of Ratner\cite{Ratner-Annals, Ratner-Acta} gives that for groups generated by one parameter unipotent subgroups any invariant probability measure on a quotient space $G/\Gamma$ has to be homogeneous. Here the rigidity is already exhibited in the action of individual one parameter subgroups of the action (another proof of this measure classification result using entropy theory was given by Margulis and Tomanov in \cite{Margulis-Tomanov}). 
Ratner used her measure classification theorem to classify orbit closures under such actions \cite{Ratner-Duke}, which enabled her to prove Raghunathan's Conjecture (this conjecture, together with a related conjecture of Dani, appeared in \cite{Dani-general-horospheric}). Several important nonhorospherical%
\footnote{The horospherical case is more elementary and can be proved e.g. using mixing of an appropriate one parameter diagonalizable flow; this is not unrelated to the phenomenon of the uniqueness of measure of maximal entropy for a one-parameter diagonalizable flow we already encountered in the context of the $\times s$ map on $\T$. We note also that the horospherical case inspired Dani and Raghunathan to make their general conjectures on unipotent orbits --- indeed, this is precisely what Dani's paper \cite{Dani-general-horospheric} is about!} cases of this conjecture were proved prior to \cite{Ratner-Duke,Ratner-Annals} by Dani and Margulis \cite{Margulis-Oppenheim-original,Dani-Margulis-primitive,Dani-Margulis-generic}, including in particular Margulis proof of the longstanding Oppenheim conjecture via the study of orbits of the group $\SO (2, 1) < \SL (3, \R)$ on $\SL (3, \R) / \SL (3, \Z)$.

There is another less important wrinkle that requires some care in formulating general conjectures regarding rigidity of higher rank abelian groups, as the following simple example illustrates: suppose for instance one would like to classify invariant measures for the action of the complex diagonal matrices on $X_\C=\SL (3, \C) / \SL (3, \mathcal{O})$ with $\mathcal{O}$ the ring of integers in an imaginary quadratic field, for instance the Gaussian integers $\Z[i]$.
$X _ \R = \SL (3, \R) / \SL (3, \Z)$, considered as a homogeneous subspace of $X _ \C$, is invariant under the real diagonal group. Let $m _ {X _ \R}$ denote the uniform measure on $X _ \R$, and set
\begin{equation*}
\mu = \fint\limits _0^ {2 \pi} \fint\limits _0^ {2 \pi} \begin{pmatrix} e ^ {i \theta _ 1}& &\\&e ^ {i \theta _ 2}&\\& & e ^ {-i (\theta _ 1+\theta _ 2)}\end{pmatrix}.m_{X _ \R} d\,\theta _ 1 d\,\theta _ 2
\end{equation*}
(here we use the symbol $\fint$ to denote integration normalized by the measure of the set we integrate on, i.e. so that $\fint \!d x =1$). The measure $\mu$ is invariant and ergodic under the action of the (complex) diagonal group in $\SL (3, \C)$, but is not homogeneous.
 
\begin{Definition}\label{definition of almost homogeneous} Let $G$ be a locally compact group, $\Gamma < G$ a closed subgroup. A measure $\mu$ on $G / \Gamma$ is said to be \emph{almost homogeneous} if there is a homogeneous measure $m _ 0$ on $G / \Gamma$ with stabilizer $H _ 0 = \stab _ G m _ 0$ and a closed subgroup $L < G$ so that $L/ (L \cap H_0)$ has finite $L$-invariant volume and
\begin{equation*}
\mu =  \fint_ {L / (L \cap H_0)} \ell . m _ 0 \,d \ell
.\end{equation*}
If the quotient $L / (L \cap H_0)$ is finite we say that $\mu$ is \emph{virtually homogeneous}.
\end{Definition}

\begin{Definition} Let $k$ be a local field (e.g. $\R$), and let $\mathbb G$ be an algebraic group defined over $k$. An element $g \in \mathbb G (k)$ said to be of \emph{class-$\mathcal{A}$} if it is diagonalizable over $k$, generates an unbounded subgroup of $\mathbb G (k)$ and moreover for any action of $\mathbb G(k)$ on a projective space $\mathbb P V (k)$ and $v \in \mathbb P V (k)$ any limit point of $\left\{ g^n .v : n \in \Z \right\}$ is $g$-invariant. An element $g \in \prod_ i \mathbb G _ i (k _ i)$ is of class-$\mathcal{A}$ if all of its components are of class-$\mathcal{A}$. 
\end{Definition}

For example, a $\R$-diagonalizable element of $\mathbb G (\R)$ with positive eigenvalues is of class-$\mathcal{A}$. Another example which works in any local field $k$ is taking an element $g \in \mathbb G (k)$ all of its eigenvalues are integer powers of some fixed $\theta \in k$ with $\absolute \theta > 1$. This latter example has been called class-$ \mathcal{A}$ by Margulis and Tomanov in \cite{Margulis-Tomanov-2}; but it is precisely the invariance property of any limit point of elements for projective actions of the underlying group that was used there, and it seems convenient to enlarge this class using this property.

\begin{Definition}\label{higher rank definition}
We say that a topological group $A$ is of \emph{higher rank} if there is a homomorphism $\Z ^ 2 \to A$ which is a proper map with respect to the discrete topology on $\Z ^2$.
\end{Definition}

General conjectures regarding rigidity for invariant measures under higher rank abelian groups were made by Furstenberg (unpublished), Katok and Spatzier \cite{Katok-Spatzier} and Margulis \cite{Margulis-conjectures}. The following is a variant of their conjectures:

\begin{Conjecture}\label{general measure conjecture} Let $\mathbb G$ be a linear algebraic group defined over $\Q$, and $S$ a finite set of places for $\Q$ containing $\infty$. Let $\mathcal{O} _ S=\Z[1/p:p\in S\setminus\infty]$ denotes the ring of $S$-integers in $\Q$, $G = \prod_ {v \in S} \mathbb G (\Q _ v)$ and $\Gamma = \mathbb G (\mathcal{O} _ S)$,\footnote{To be more precise: we fix a realization of $\mathbb G$ as a $\Q$-subgroup of $\SL (d)$ for some $d$, and set $\Gamma = \mathbb G (\Q) \cap \SL (d, \mathcal{O} _ S)$.} diagonally embedded in $G$.
Let $A<G$ be a closed subgroup consisting of elements of class-$\mathcal A$ of higher rank and let $\mu$ be an $A$-invariant and ergodic probability measure on $G / \Gamma$.
Then either $\mu$ is virtually homogeneous or there is a $\Q$-subgroup $\mathbb L \leq \mathbb G$ and a proper normal $\Q$-subgroup $\mathbb H \lhd \mathbb L$ so that, if $H = \prod_ {v \in S} \mathbb H(\Q _ d)$ and $L = \prod_ {v \in S} \mathbb L(\Q _ d)$, then
\begin{enumerate}
\item $A \cap L$ has finite index in $A$,
\item there is some $g \in G$ so that $\mu (g.[L]_\Gamma)>0$ (with $[\bullet] _ \Gamma$ denoting the image under the projection $G \to G / \Gamma$)\footnote{We will also use the notation $[g]$ for $[g] _ \Gamma$ when $\Gamma$ is understood.}
\item the image of $A \cap L$ in $L / H$ is \textbf{not} of higher rank.
\end{enumerate}
\end{Conjecture}

Unlike the case of unipotent flows, where the classification of invariant measures and orbit closures go hand-in-hand and are very closely analogous, for diagonal flows the problem of classifying invariant measures seems, in general, better behaved than understanding orbit closures. This is somewhat surprising, as in the $\times a, \times b$-system considered by Furstenberg, with $a$ and $b$ multiplicatively independent integers, a complete orbit closure classification was obtained by Furstenberg already in 1967, whereas the measure classification question (without a positive entropy assumption) is Conjecture~\ref{Furstenberg's conjecture} --- a notoriously hard open problem. However, already in the slight generalization of considering $\times a, \times b$ for $a$ and $b$ multiplicatively independent (rational) integers on $\C/\Z[i]$ not much is known (see \S\ref{torus orbits section} for more details).

While this is not immediately clear from the formulation, Conjecture~\ref{Furstenberg's conjecture} is essentially a special case of Conjecture~\ref{general measure conjecture}. For simplicity, assume that $a$ and $b$ are distinct primes (the modification for general multiplicatively independent $a$ and $b$ is left to the imagination of the reader). Let $\mathbb G = \left\{ \begin{pmatrix}*&*\\ 0& 1 \end{pmatrix} \right\}$, i.e. the semi-direct product of the multiplicative group $\mathbf G_m$ with the additive group $\mathbf G_a$, and take $S= \left\{ \infty, a, b \right\}$. Let $G = \prod_ {v \in S} \mathbb G (\Q _ v)$, $\Gamma = \mathbb G (\Z [1/ab])$ and $A<\prod_ {v \in S} \mathbf G_m (\Q _ v)<G$ the group $\left\{ a^kb^l:k,l \in \Z \right\}$ diagonally embedded in $\prod_ {v \in S} \mathbf G_m (\Q _ v)$. For any $y \in \R ^{\times}$ we have that
\begin{equation*}
Y_y := \left\{ \left [\begin{pmatrix}y&x_\infty\\ 0& 1 \end{pmatrix} , \begin{pmatrix}y_a&x_a\\ 0& 1 \end{pmatrix}, \begin{pmatrix}y_b&x_b\\ 0& 1 \end{pmatrix}  \right] _ \Gamma \,\middle| \begin{array}{l}\text {$y_v \in \Z _ v ^{\times}$ for $v=a,b$}\\ \text{$x _ v \in \Q _ v$ for $v =a,b, \infty$}\end{array} \right\}
\end{equation*}
is a compact $A$-invariant subset of $G / \Gamma$, hence any $A$-invariant and ergodic probability measures $\mu$ on $G / \Gamma$ is supported on a single $Y _ y$. Without loss of generality we can assume it is supported on $Y_1$. Let $\pi$ be the projection of $Y_1$ to $X_{a,b}=\R \times \Q _ a \times \Q _ b / \Z [1/ab]$ given by
\begin{equation*}
\left [\begin{pmatrix}1&x_\infty\\ 0& 1 \end{pmatrix} , \begin{pmatrix}y_a&x_a\\ 0& 1 \end{pmatrix}, \begin{pmatrix}y_b&x_b\\ 0& 1 \end{pmatrix} \right] _ \Gamma \mapsto \bigl[x_\infty,x_a,x_b\bigr]_{\Z [1/ab]}
.\end{equation*}
For any $A$-invariant and ergodic probability measure $\mu$ supported on $Y _ 1$, the measure $\pi _*\mu$ is a $\times a, \times b$-invariant and ergodic probability measures on $X_{a,b}$. Conversely, since the fibers of the map $\pi: Y _ 1 \to X_{a,b}$ are compact, any $\times a, \times b$-invariant and ergodic probability measures on $X_{a,b}$ can be lifted to a $A$-invariant and ergodic measure on $Y_1$.

In this survey we focus on \emph{$S$-arithmetic quotients}: quotients of a finite index subgroup $G$ of the $\Q_S$ points $\mathbb G(\Q_S)$ of a $\Q$-group $\mathbb G$ by a subgroup $\Gamma$ commensurable to $\mathbb G(\mathcal O_S)$.  By restriction of scalars, this implicitly also includes the case of algebraic groups defined over any number field, but because of issues related to those pointed out above for $\SL (3, \C) / \SL (\Z [i])$, it is more convenient to work with the smaller field $\Q$. 

\begin{Definition} An $S$-arithmetic quotient $G / \Gamma$ is \emph{saturated by unipotents} if it has finite volume and the group generated by one parameter unipotent subgroups of $G$ acts ergodically on $G /  \Gamma$ (with respect to the Haar measure on $G / \Gamma$).\label{saturated by unipotents}
\end{Definition}

\noindent
When working with real algebraic groups $\mathbb G (\R)$ where $\mathbb G(\C)$ is generated by unipotents (equivalently, the radical of $\mathbb G$ is equal to the unipotent radical of $\mathbb G$), a quotient $G/\Gamma$ satisfies the  saturated by unipotents property if and only if  it is connected in the Hausdorff topology (cf. \cite[Ch.~II]{Margulis-book}).

\medskip

A very interesting and active direction we do not cover in this survey is actions on quotients of algebraic groups defined over global fields of positive characteristic. The key feature here is that there is no analogue to ``$ \Q$'': there is no minimal global field. This type of issue makes analyzing even the analogue of the $\times a, \times b$-system in positive characteristic quite intricate (cf.\ \cite[Construction 5.2]{Kitchens-Schmidt} and \cite{Einsiedler-bad-measures}). For quotients of semisimple groups the situation is better, and a measure classification theorem for positive entropy measures analogous to what Einsiedler, Katok and the author \cite{Einsiedler-Katok-Lindenstrauss} proved for $\Q$  has been proved by Einsiedler, Mohammadi, and the author in \cite{Einsiedler-Lindenstrauss-Mohammadi}. However, even in this case, it is far from clear to which extent one should expect an analogue to Conjecture~\ref{general measure conjecture}; in this context we mention the paper \cite{Adiceam-Nesharim-Lunnon} by Adiceam, Nesharim and Lunnon where a very interesting example is constructed.

\subsection*{Acknowledgements} This paper is dedicated with admiration to Gregory Margulis whose deep and profound work has been, and continues to be, an inspiration to me ever since I started getting interested in homogeneous dynamics. Indeed, Margulis' deep work has been a big part of what drew me to the subject to begin with. On a more personal level, I would like to thank him for his kindness and generosity over the years. 

Many of the results I describe in this work are joint with Manfred Einsiedler; it is a pleasure to express my gratitude to him for this collaboration; I would also like to thank him for comments on earlier versions of this survey. I also thank Ilya Khayutin for helpful comments and corrections.
Finally, I would like to thank the editors of this volume for inviting me to contribute to it and for their patience.

\section{Measure rigidity of higher rank diagonal actions}\label{measure rigidity section}
While Conjectures \ref{Furstenberg's conjecture} and \ref{general measure conjecture} are still wide open, significant progress was obtained regarding classifying invariant measures under a positive entropy condition. In sections \S\ref{sec: applications 1}--\ref{sec: applications 3} we survey some applications of these results. Typically we are given a sequence of $A$-invariant probability measures $\mu _ i$ on $G / \Gamma$, and would like to understand what are the weak$^ {*}$ limit points of the sequence $\mu _ i$. Suppose $\mu$ is such a limit. A priori it seems very difficult to control any kind of ergodicity or mixing condition for the limiting measure. On the other hand, entropy is fairly well behaved with respect to weak$^ {*}$ limits. For example, for an $A$-invariant measure $\nu$, let $h(\nu, a)$ denote the ergodic theoretic (a.k.a Komogorov-Sinai or metric) entropy of  $\nu$ with respect to the action $[g] _ \Gamma \mapsto a.[g] _ \Gamma$. Then if $G / \Gamma$ is compact, if $\mu _ i \to \mu$ weak$^*$ then for any $a \in A$
\begin{equation*}
h (\mu, a) \geq \limsup\nolimits_ i h (\mu _ i ,a)
\end{equation*}
(cf.\ e.g.~\cite[\S9]{Einsiedler-Katok-Lindenstrauss}).
This actually also holds if $G / \Gamma$ is not compact \textbf{assuming} $\mu$ is a probability measure.

Rudolph's theorem \cite{Rudolph-2-and-3} discussed above (p.~\pageref{Rudolph theorem mention}) regarding $S_{a,b} = \left\{ a ^ k b ^ l: k, l \in \N \right\}$-invariant and ergodic measure $\mu$ on $\T$ has been a prototype for many subsequent theorems. We remark that a simple yet important lemma in Rudolph's proof implies that if $h (\mu, s) > 0$ for one $s \in S_{a,b}$, then $h (\mu, s) > 0$ for all $s \in S_{a,b}$. Katok and Spatzier (\cite{Katok-Spatzier, Katok-Spatzier-corrections}, cf.\ also \cite{Kalinin-Katok-Seattle} by Kalinin and Katok) pioneered the study following Rudolph of higher rank abelian actions by automorphisms on $\T ^ d$ and by translations on quotients $G / \Gamma$ (note that similarly to what we have seen for the $\times a, \times b$-case, the former can be viewed as a special case of the latter where the group $\mathbb G$ is a semi-direct product of a torus and an abelian additive group). In some cases, Katok and Spatzier were able to obtain a full analogue of Rudolph's theorem, but in most cases (e.g. for $\Z^k$ actions on $G/\Gamma$ with $G$ semisimple) an additional ergodicity condition is needed, a condition that unfortunately is not stable under weak$^*$ limits.

\subsection{Rigidity of joinings}
Arguably the most complete result regarding the classification of higher rank abelian actions on arithmetic quotients does not explicitly mention entropy, though entropy plays an important role in the proof. In the same paper \cite{Furstenberg-disjointness-1967} in which Furstenberg proved Theorem~\ref{Furstenberg's theorem}, thereby introducing higher rank rigidity from the dynamical perspective, Furstenberg also introduced joinings as a key tool in the study of measure preserving and topological dynamical systems. Suppose $H$ is a topological group acting in a measure preserving way on two probability measures spaces $(X, \mu)$ and $(X ', \mu ')$. Then $H$ also acts on the product space $X \times X'$ by setting $h.(x,x')=(h.x,h.x')$. A \emph{joining} of $(X, \mu, H)$ and $(X ', \mu ', H )$ is an $H$-invariant probability measure on $X \times X '$ that projects to the measure $\mu$ and $\mu '$ on $X$ and $X '$ respectively. There is always at least one joining between any two such actions, namely the product measure $\mu \times \mu '$. Existence of other joinings can be interpreted as evidence of some communality between $(X, \mu, H)$ and $(X ', \mu ', H )$; an extreme form of this would be if these two measure preserving $H$-actions would be isomorphic (as $H$-actions, i.e. there is a measure preserving 1-1 and onto map $\phi$ between subsets of full measure of $X$ and $X '$ commuting with the $H$-action), in which case the push forward under $(\mathrm{id}, \phi)$ of $\mu$ would be a nontrivial joining supported on the graph of~$\phi$.

The following general joining classification theorem is the main result of \cite{Einsiedler-Lindenstrauss-joinings-2} by Einsiedler and the author:

\begin{Theorem}[\cite{Einsiedler-Lindenstrauss-joinings-2}] 
\label{higher rank joining}
Let $r,d\geq 2$ and let $\GG_1,\ldots,\GG_r$ be
semisimple algebraic groups defined over~$\Q$ that are~$\Q$-almost simple, $\GG=\prod_{i=1}^r\GG_i$, and $S$ be a finite set of places of $\Q$.
Let $X_i = \Gamma _ i \backslash G _ i$ be $S$-arithmetic quotients\footnote{In particular, by our definition of $S$-arithmetic quotients $G_i$ has finite index in $\GG_i(\Q_S)$.}  saturated by unipotents for $G_i \leq \GG_i(\Q_S)$  and let $X= \prod_ i X_i$. Let $a_i : \Z ^ d \to G _ i$ be proper homomorphisms so that $a=(a_1, \dots, a_r) : \Z ^ d \to G= \prod_ i G_i$ is of class-$\mathcal{A}$, and set $A = a(\Z ^ d)$. Suppose $\mu$ is an $A$-invariant and ergodic joining of the actions of $A_i = a_i(\Z^d)$ on $X_i$ equipped with the Haar measure $m_{X_i}$. Then $\mu$ is homogeneous.
\end{Theorem}
\noindent
In fact, \cite{Einsiedler-Lindenstrauss-joinings-2} gives slightly more precise information, in that $\mu$ is not just homogeneous but Haar measure on a finite index subgroup of the $S$-adic points of an algebraic group defined over $\Q$. Such a measure would be said to be \emph{an algebraic measure defined over $\Q$}. This joining classification theorem can be extended to perfect groups. Recall that an algebraic group $\mathbf G$ is said to be perfect if $\mathbf G = [\mathbf G, \mathbf G]$.

\begin{Theorem}[\cite{Einsiedler-Lindenstrauss-joinings-2}] \label{thm:perfect joining}
Let $r,d\geq 2$ and let $\GG_1,\ldots,\GG_r$ be perfect algebraic groups defined over~$\Q$, $\GG=\prod_{i=1}^r\GG_i$, and $S$ be a finite set of places of $\Q$.
Let $X_i = \Gamma _ i \backslash G _ i$ be $S$-arithmetic quotients for $G_i \leq \GG_i(\Q_S)$ saturated by unipotents and let $X= \prod_ i X_i$. Let $a_i : \Z ^ d \to G _ i$ be homomorphisms so that $a=(a_1, \dots, a_r) : \Z ^ d \to G= \prod_ i G_i$ is of class-$\mathcal{A}$, and such that the projection of $a_i$ to every $\Q$ almost simple factor of $\GG_i(\Q_S)$ is proper.
Suppose $\mu$ is an $A$-invariant and ergodic joining of the action of $A_i = a_i(\Z^d)$ on $X_i$ equipped with the Haar measure $m_{X_i}$. Then $\mu$ is homogeneous, indeed an algebraic measure defined over $\Q$.
\end{Theorem}

We remark that for the action of a one parameter unipotent group on quotients of $\SL(2,\R)$ by lattices Ratner established a joining classification theorem in \cite{Ratner-products}. A general joining classification result for actions of unipotent groups was given by Ratner in \cite{Ratner-Acta} as a by product of her techniques to classify all invariant measures\footnote{Indeed, the joining classification follows directly from the measure classification theorem of Ratner in \cite{Ratner-Annals}, but in \cite{Ratner-Acta} (which is part of the sequence of papers establishing the results \cite{Ratner-Annals}) this result is already noted.}.

The restriction to perfect groups in Theorem~\ref{thm:perfect joining} is important. If $\alpha$ is a (faithful) $\Z ^ k$-action on a torus $\T ^ d$ by automorphisms, or more generally a $\Z ^ k$-action on the solenoid $\T_S ^ d = \left (\prod_ {v \in S} \Q _ v \right) ^ d$ (a prime example of the later being the action generated by the $\times a$ and $\times b$ maps on $\T _ S$ with $S$ containing $\infty$ as well as all prime factors of $a b$) for $k \geq 2$ then any hypothetical nonatomic $\alpha (\Z ^ d)$-invariant and ergodic invariant measure on $\T _ S^d$ of zero entropy would give rise to a nontrivial, nonhomogeneous, self joining of $(\T _ S^d, m_{\T _ S^d}, \alpha (\Z ^ k))$ given by the push forward of the measure $m_{\T _ S^d} \times \mu$ using the map $(x,y) \mapsto (x, x+y)$ from $\T _ S^{2d} \to \T _ S^{2d}$.
This simple example shows that classifying self joinings of such $\Z ^ k$-actions is (at least) as hard as Conjecture~\ref{Furstenberg's conjecture}. However, one can classify joinings between such $\Z ^ d$-actions up to zero entropy quotients (\cite{Kalinin-Katok,Kalinin-Spatzier,Einsiedler-Lindenstrauss-ERA}).

\subsection{Some measure classification theorems for $S$-arithmetic quotients}\label{higher rank measure rigidity}

Joinings between higher rank abelian actions have positive entropy coming from the factors being homogenous, but in fact being a joining imposes additional restrictions on leafwise measures that are very useful for the analysis. If one wants a measure classification of positive entropy measures, some additional conditions are needed.

One condition which gives rise to a clean statement is when the acting group is a maximal split torus, or more generally satisfies the following condition:

\begin{Definition}
Let $\mathbb G$ be an algebraic group defined over $\Q$, $S$ a finite set of places. A subgroup $A < \mathbb G (\Q _ S)$ will be said to be a \emph{partially maximal $\Q _ S$-split torus} if there is for each $s \in S$ a (possibly trivial) algebraic normal subgroup\footnote {To be precise, $\mathbf H _ s$ is a group of $\Q _ s$-points of a $\Q _ s$ group, however (in contrast to the global field case) when considering an algebraic group over local field will not make the distinction between the abstract algebraic group and the groups of points.} $\mathbf H_s \lhd \mathbb G (\Q _ s)$ so that $(A \cap \mathbf H _ s)$ is a maximal $\Q _ s$-split torus in $\mathbf H _ s$, and $A = \prod_ {s \in S '} (A \cap \mathbf H _ s)$.
\end{Definition}

\begin{Theorem}[Einsiedler and L.~\cite{Einsiedler-Lindenstrauss-split}]\label{split theorem simple}
Let $\mathbb G$ be a $\Q$-almost simple algebraic group, $S$ a finite set of places, and $G/\Gamma$ an $S$-arithmetic quotient for $\mathbb G$ saturated by unipotents in the sense of Definition~\pageref{saturated by unipotents}. Let $A$ be a higher rank, partially maximal $\Q _ S$-split torus.
Let $\mu$ be a $A$-invariant and ergodic measure on $G / \Gamma$, and assume that: 
\begin{enumerate}
\item $\mu(g [\mathbb L (\Q _ S) \cap G]_\Gamma)=0$ for every proper reductive subgroup $\mathbb L < \mathbb G$ and $ g \in G$;
\item $h (\mu, a) > 0$ for some $a \in A$.
\end{enumerate}
Then $\mu$ is the uniform measure on $G / \Gamma$.
\end{Theorem}

\noindent
Using Theorem~\ref{split theorem simple} a decomposition theorem can be proved for measures on $S$-arithmetic quotients corresponding to a semisimple $\Q$-group as a product of four pieces, that may well be trivial:

\begin{Theorem}[Einsiedler and L.~\cite{Einsiedler-Lindenstrauss-split}]\label{split theorem semisimple}
Let $\mathbb G$ be a semisimple algebraic group defined over $\Q$, $S$ a finite set of places, $G/\Gamma$ an $S$-arithmetic quotient for $\mathbb G$, and $A$ a partially maximal $\Q _ S$-split torus.
Let $\mu$ be a $A$-invariant and ergodic measure on $G / \Gamma$. Then there is a finite index subgroup $A ' < A$ and a probability measures $\mu '$ so that $\mu = \frac {1 }{ \absolute {A / A '}} \sum_ {a \in A / A '} a \mu ' $ and so that $\mu '$ can be decomposed as follows. For $i \in \left\{ 1,2,3 \right\}$ there is a semisimple $\Q$-subgroup $\mathbb L _ i \leq \mathbb{G}$ and an anisotropic $\Q$-torus $\mathbb L _ 0 < \mathbb{G}$ so that $\iota: (l_0, \dots, l_3) \mapsto l_0 \cdot \dots l_3$ gives a finite to one map $ \prod_ {i = 0} ^ 3 \mathbb L _ i (\Q _ S) \to \mathbb G (\Q _ S)$ so that $\mu ' = \iota _ {*} (\mu_0 \times \dots \times \mu_3)$ with each $\mu _ i$ an $A ' \cap \mathbb L _ i (\Q _ S)$-invariant and ergodic probability measure on $(\mathbb L _ i (\Q _ S) \cap G) / (\mathbb L _ i (\Q _ S) \cap \Gamma)$, $A' = \prod_ {i = 0} ^ 3 (A \cap \mathbb L _ i (\Q _ S))$ and
\begin{enumerate}
\item $\mu _ 1$ is the uniform measure on $L / (\mathbb L _ 1 (\Q _ S) \cap \Gamma)$ with $L \leq \mathbb L _ 1 (\Q _ S)$ a finite index subgroup,
\item $\mu _ 2$ satisfies that $h (\mu _ 2, a) = 0$ for every $a \in A \cap \mathbb L _ 2 (\Q _ S)$,
\item $\mathbb L _ 3 (\Q _ S)$ is an almost direct product of $\Q$-almost simple groups $\mathbb L _ {3,i} (\Q _ S)$ so that for all $i$ the group $A \cap \mathbb L _ {3,i} (\Q _ S)$ is not of higher rank.
\end{enumerate}
\end{Theorem}

The special cases of Theorems~\ref{split theorem simple} and \ref{split theorem semisimple} for $G/\Gamma$ a quotient of $\prod_ {i=1}^k \SL (2, \Q_{v_i})$ by an irreducible lattice 
 (with $v_i \in \{\text{primes or $\infty$}\}$) or $G/\Gamma= \SL (n, \R)/ \SL (n, \Z) $ 
 were proven earlier by the author \cite{Lindenstrauss-quantum} and Einsiedler, Katok and the author~\cite{Einsiedler-Katok-Lindenstrauss} respectively.
 
 Ideally, one would like to obtain a measure classification results for measures invariant under a higher rank diagonalizable group in the more general context of Conjecture~\ref{general measure conjecture}. This is the subject of ongoing work; in particular, in joint work with Einsiedler we have the following:

\begin{Theorem}[Einsiedler and L.~\cite{Einsiedler-Lindenstrauss-general-SL2}]\label{sl2-thm-final}
Let~$\GG$ be an algebraic group over~$\Q$ that is~$\Q$-almost simple
and a form of~$\operatorname{SL}_2^k$ or of~$\PGL_2^k$ with~$k\geq 1$, $S$ a finite set of places, and $G/\Gamma$ an $S$-arithmetic quotient for $\mathbb G$.
Let~$A<G$ be a closed abelian subgroup of class-$\cA$ and of higher rank.  Let~$\mu$
be an~$A$-invariant and ergodic probability measure on~$X=\Gamma\backslash G$ such that~$h_\mu(a)>0$
for some~$a\in A$.  Then one of the following holds:
\begin{itemize}
\item \textup{\bf (Algebraic)} the measure~$\mu$ is homogenous.
\item \textup{\bf (Solvable)} the space~$X$ is non-compact.
There exists a nontrivial unipotent subgroup~$L$
such that~$\mu$ is invariant under~$L$. The measure $\mu$ is
supported on a compact~$A$-invariant orbit~$x_0M\cong\Lambda_{M,x_0}\backslash M$,
where~$M<G$ is a solvable subgroup and~$\Lambda_{M,x_0}=\{m\in M\mid m.x_0=x_0\}$
is the stabilizer of~$x_0$ in~$M$. The lattice $\Lambda_{x_0,M}$ in $M$ intersects
the normal subgroup~$L\lhd M$ in a uniform lattice
and if~$\pi:M\rightarrow M/L$ denotes the natural projection map, then
the image of~$\mu$ under the
induced map~$\Lambda_{M,x_0}\backslash M\rightarrow \pi(\Lambda_{M,x_0})\backslash (M/L)$
has zero entropy for the action of~$A$.
\end{itemize}
\end{Theorem}

 \subsection{A rigidity theorem for measures invariant under a 1-parameter diagonal group with an additional recurrence assumption}
 
 For the application of measure rigidity to quantum unique ergodicity, a variant of the above results was essential, where the assumption of invariance under a higher rank group was relaxed.

\begin{Definition} Let $H$ be a locally compact group acting on a standard Borel space $(X, \mathcal B) $. We say that a measure $\mu$ on $X$ is $H$-recurrent\footnote{An alternative terminology often used in this context is $H$-conservative; we prefer $H$-recurrent as it seems to us more self-explanatory.} if for every $B \subset X$ with $\mu (B) > 0$ and any compact subset $F \subset H$, for $\mu$-a.e. $x \in X$ there is an $h \in H \setminus F$ with $h.x \in B$.
\end{Definition}

\noindent
We stress that no assumption is made regarding $H$-invariance of $\mu$ or even the measure class of $\mu$.

\begin{Theorem} [L.~\cite{Lindenstrauss-quantum}]\label{one parameter plus recurrence rigidity theorem}
Let $G = \prod_ {i=1}^r \SL (2, \Q _ {v_i})$ with $v_i \in \left\{ \infty, \text{primes} \right\}$ and $r \geq 2$, and let $\Gamma < G$ be an irreducible lattice. Let $A<\SL (2, \Q _ {v_1})$ be a 1-parameter diagonal group, $a \in A$ generating an unbounded subgroup of $A$, and  $H = \prod_ {i=2}^r \SL (2, \Q _ {v_i})$. Suppose $\mu$ is a $A$-invariant, $H$-recurrent, and that for a.e.~$A$-ergodic component $\mu _ \xi$ of $\mu$ the entropy $h (\mu _ \xi, a)>0$. Then $\mu$ is the uniform measure on $G / \Gamma$.
\end{Theorem}

We note that using recurrence as a substitute for invariance under a higher rank group was motivated by Host's proof of Rudolph's Theorem in \cite{Host-normal-numbers}.

\section{Orbit closures - many questions, a few answers} \label{sec:orbit closures}

\subsection{Prologue - orbit closures and equidistribution for unipotent flows}\label{unipotent section}
For unipotent flows, there is a very close relationship between behavior of individual orbits and the ergodic invariant measures. This correspondence was used by Ratner  \cite{Ratner-Duke} to prove the Raghunathan conjecture: 

\begin{Theorem}[Ratner \cite{Ratner-Duke}]
Let $U$ be a connected unipotent subgroup of real algebraic group $G$ and $\Gamma < G$ a lattice. Then for any $x \in G / \Gamma$ there is a closed subgroup $U \leq L \leq G$ so that  $\overline {U . x}=L.x$, with $L . x$ a periodic orbit (i.e., $\stab _ L (x)$ is a lattice in $L$). Moreover, $U$ acts ergodically on $L / \stab  _ L (x)$.
\end{Theorem}

\noindent
In particular, the orbit closure of every $U$-orbit $U.x$ is the support of a $U$-invariant and ergodic measure on $G / \Gamma$.

One key ingredient used to prove this surprisingly tight correspondence is a nondivergence estimate for unipotent flows developed by Dani and Margulis \cite{Margulis-nondivergence,Dani-nondivergence,Dani-Margulis-nondivergence}.
In addition to establishing nondivergence of the $U$-trajectory, needed in order to obtain from an orbit some limiting probability measure which can be analyzed, to deduce the Raghunathan's Conjecture from the measure classification theorem one needs to establish that a trajectory of a point $x \in G / \Gamma$ does not spend a lot of time close to a ``tube'' corresponding to shifts of a given periodic $L$ orbit, unless $x$ itself is in this tube. A flexible way to establish such estimates, known as the \emph{Linearization Method}, was developed by Dani and Margulis \cite{Dani-Margulis-linearization}; while \cite{Dani-Margulis-linearization} uses Ratner's measure classification theorem, the technique itself was developed earlier by Dani and Margulis (with closely related works by Shah) in order to prove some cases of Raghunathan's conjecture by purely topological means (see e.g.~\cite{Dani-Margulis-generic}); an alternative approach to linearization was used by Ratner in her proof of Raghunathan Conjecture. We mention that a stronger (more) explicitly effective version of the Dani-Margulis Linearization Method was given recently by Margulis, Mohammadi, Shah and the author \cite{Lindenstrauss-Margulis-Mohammadi-Shah}.

We also recall the following theorem of Mozes and Shah that relies on Ratner's measure classification theorem and the linearization method:

\begin{Theorem} [Mozes and Shah \cite{Mozes-Shah}]\label{Mozes-Shah theorem} Let $G$ be a linear algebraic group over $\R$, $\Gamma < G$ a lattice, and let $\mu _ i$ be a sequence of probability measures on $G / \Gamma$ and $u ^ {(i )} _ t$ a sequence of one parameter unipotent subgroups of $G$ so that for every $i$ the measure $\mu _ i$ is $u ^ {(i )} _ t$-invariant and ergodic. Suppose $\mu _ i$ converges in the weak$^{*}$-topology to a probability measure $\mu$. Then $\mu$ is homogeneous, and moreover there are $g _ i \to e$ so that $g _ i \supp (\mu _ i) \subset \supp (\mu)$ for $i$ large enough.
\end{Theorem}

\noindent
This theorem was extended to the $S$-arithmetic setting by Gorodnik and Oh \cite{Gorodnik-Oh-adelic-periods}.

\subsection{Orbits closures for higher rank diagonalizable group in a torus}\label{torus orbits section}

Actions of one parameter diagonal groups display no rigidity, and most questions about behavior of individual orbits for one parameter diagonal groups seem to be hopelessly difficult. For instance, it is a well-known open problem whether $\sqrt[3]{2}$ (or indeed any other irrational algebraic number of degree $\geq 3$) has a bounded continued fraction expansion, which is completely equivalent to the question whether the half-orbit
\begin{equation*}
\left\{ \begin{pmatrix} e ^ t& \\& e ^ {- t} \end{pmatrix} \left [\begin{pmatrix}  1&\sqrt[3]{2} \\0& 1 \end{pmatrix} \right] : t \geq 0 \right\}
\end{equation*}
is bounded in $G / \Gamma = \SL (2, \R) / \SL (2, \Z)$.

One could hope that the situation would be better for actions of higher rank diagonal groups, which do have some rigidity, and to a certain extent this is true. However, any hope of obtaining as good an understanding of orbits of higher rank diagonal groups as we have for unipotent flows is doomed to failure, in large parts stemming from the fact that the connection between individual orbits and invariant measures for diagonal flows is much weaker.

To illustrate this point, consider first Furstenberg's Theorem \ref{Furstenberg's theorem}. This theorem gives a complete classification of orbit closures for the action of $S _ {a, b} = \left\{ a ^ n b ^ m: n, m \in \N \right\}$ on $\T = \R/ \Z$: Either a finite orbit on which $S _ {a, b}$ acts transitively, or $\T$. This was significantly extended by Berend, who showed the following

\begin{Theorem} [Berend, \cite{Berend-invariant-tori, Berend-invariant-groups}]\label{Berend's theorem}
Let $K$ be a number field, $S$ a finite set of places including all infinite places, and $\mathcal{O} _ S$ the ring of $S$-integers, i.e. the ring of $k \in K$ satisfying that $\absolute k \leq 1$ for any place $v \not\in S$ of $K$. Let $\Sigma$ a higher rank subgroup of $\mathcal{O} _ S ^{*}$. Assume that
\begin{enumerate}
\item no finite index subgroup of $\Sigma$ is contained in a proper subfield of $K$,
\item for every $v \in S$ there is some $a \in \Sigma$ with $\absolute a _ v >1$
.\end{enumerate}
Then any $\Sigma$-invariant closed subset of $X=\prod_ {v \in S} K _ v / \mathcal{O} _ S$ is either finite or $X$ itself.
\end{Theorem}

For $K = \Q$ and $\Sigma \subset \N $ (not including 0!) this reduces easily to Furstenberg theorem.
The proofs of Furstenberg and Berend (as well as a simple proof of Furstenberg's Theorem by Boshernitzan \cite{Boshernitzan-elementary}) are purely topological, but one can deduce Theorem~\ref{Furstenberg's theorem} and Theorem~\ref{Berend's theorem} from Rudolph's theorem and its analogue to solenoids \cite{Einsiedler-Lindenstrauss-ERA} by Einsiedler and the author respectively by establishing that any infinite closed $\Sigma$-invariant subset  $Y \subset X=\prod_ {v \in S} K _ v / \mathcal{O} _ S$ has to support a $\Sigma$-invariant measure of \emph{positive entropy}. The reason this can be shown is that it is not hard to show that if $Y$ is such a closed, infinite, invariant subset $Y - Y = X$. This approach was used by Bourgain, Michel, Venkatesh and the author to give a quantitative version of Furstenberg theorem in \cite{Bourgain-Lindenstrauss-Michel-Venkatesh}.

Both conditions in Theorem~\ref{Berend's theorem} are needed in order to ensure that any closed invariant subset is either finite or $X$. However, dropping assumption (2) does not dramatically change the situation: if there is some $v \in S$ so that for every $a \in \Sigma$ we have that $\absolute {a} _ v = 1$ then there would certainly be other possible orbits closures, e.g.\ a $\Sigma$-invariant subset of $X$ supported on the $K _ v$-orbit of the origin on which $\Sigma$ acts by generalized rotations. However, with the minor necessary changes needed to accomodate such obvious examples of orbit closures, the above classification holds also without assumption (2). This was shown by Wang \cite{Wang-non-hyperbolic} for the case of $X$ being a torus (and essentially the same proof also works for the more general class of $X$ considered in Theorem~\ref{Berend's theorem}); we also mention that a very interesting combinatorial applications for the case of $K= \Q (i)$ was given by Manners in \cite{Manners-pyjama} (who gave an independent treatment of the relevant orbit closure classification theorem). We also note that the approach outlined in the previous paragraph to proving Theorem~\ref{Berend's theorem} using the measure classification result of Einsiedler and the author \cite{Einsiedler-Lindenstrauss-ERA} works just as well in the case where assumption (2) does not necessarily hold\footnote{The paper \cite{Einsiedler-Lindenstrauss-ERA} gives a full treatment of a measure classification theorem assuming positive entropy for irreducible actions, which is the case relevant here, as well as announces results for more general cases with some hints regarding proofs.}.

If one instead weakens the conditions of Theorem~\ref{Berend's theorem} by eliminating the irreducibility assumption (1) (even keeping assumption~(2)) we already enter the realm of conjectures, where surprisingly difficult questions loom. For instance, suppose $\Sigma$ is contained in a subfield $L < K$ with $[K: L] = 2$, but that no finite index subset of $\Sigma$ is contained in a proper subfield of $L$.
Suppose even that $\Sigma$ is the full group of units of the $S$-integers of $L$ (or more precisely, the $S_L$-units of $L$, with $S_L$ the set of places of $L$ corresponding to the those in $S$), and $S$ consists only of all infinite places of $K$ (so that $X$ is a torus\footnote{The assumption that $\Sigma$ is the full group of units of the $S$-integers of $L$ is a significant assumption; the assumption that $X$ is a torus, i.e.~$S$ consists only of all infinite places of $K$, can be removed.}). Then if the rank of $\Sigma$ is $\geq 3$ Wang and the author \cite{Lindenstrauss-Wang} proved that any orbit closure is (at most) a finite union of cosets of closed (additive) subgroups of $X$. Surprisingly, this statement is false for $\Sigma$ of rank 2! We make however the following conjecture:

\begin{Conjecture}\label{torus conjecture}
Let $K$ be a number field, $S$ a finite set of places including all infinite places, and $\mathcal{O} _ S$ the ring of $S$-integers, i.e. the ring of $k \in K$ satisfying that $\absolute k \leq 1$ for any place $v \not\in S$ of $K$. Let $L<K$, $S_L$ the set of valuations of $L$ corresponding to places in $S$, and $\Sigma$ a higher rank subgroup of the group of $S_L$-units of $L$; assume moreover that no finite index subgroup of $\Sigma$ is contained in a proper subfield of~$L$. Let $X=\prod_ {v \in S} K _ v / \mathcal{O} _ S$ and
let $Y= \overline {\Sigma .x}$ for $x \in X$. Then either $Y=X$ or there exists a finite collection $X_i$ of closed proper subgroups of $X$ and torsion points $p_i \in X$ so that $Y \subset \Sigma . x \cup \bigcup_i (X _ i + p _ i)$.
\end{Conjecture}

We remark that (at least when $[K:L] = 2$) one can give a complete classification of the \emph{support} of $\Sigma$-ergodic and invariant measures, and (at least to us) it seems that the key difficulty in proving Conjecture~\ref{torus conjecture} is the weak correspondence between invariant measures and individual orbits in the diagonalizable case, in sharp contrast to \S\ref{unipotent section}.

Conjecture~\ref{torus conjecture} is somewhat close in spirit to a recent result of Peterzil and Strachenko \cite{Peterzil-Starchenko-torus} (which they extended later to nilmanifolds) that proves a similar structure for the image of a definable subset of $\R ^ d$ with respect to an $o$-minimal structure in $\T ^ d = \R ^ d / \Z ^ d$.

\subsection{Orbit closures and limits of periodic measures for actions of higher rank diagonalizable groups on quotients of semisimple groups}

We already mentioned in the introduction the important conjecture of Cassels and Swinnerton-Dyer regarding orbit closures of the full diagonal group $A$ in the homogeneous space $\SL (d, \R) / \SL (d, \Z)$.\footnote{We implicitly identify between $\SL (d, \R) / \SL (d, \Z)$ and $\PGL (d, \R) / \PGL (d, \Z)$; while the underlying algebraic groups are different, the quotients are isomorphic.} For the convenience of the reader, we recall it here:

\begin{ConjectureA}[\cite{Cassels-Swinnerton-Dyer,Margulis-Oppenheim-conjecture}]
Let $d \geq 3$. Any orbit of the diagonal group $A$ in $\PGL (d, \R) / \PGL (d, \Z)$ is either unbounded or periodic.
\end{ConjectureA}

One would liked to say at least conjecturally something stronger about the orbit closure of an orbit $A . x$ for $A . x$ nonperiodic. For instance, in the same paper Cassels and Swinnerton-Dyer give a conjecture which can be phrased as saying that any orbit of $\SO (2,1)$ on $\SL (3, \R) / \SL (3, \Z)$ is either periodic or unbounded, a conjecture that is a special case of Raghunathan's conjecture and was proved by Margulis in the mid 1980s \cite{Margulis-Oppenheim-CR, Margulis-Oppenheim-original}. As we saw in \S\ref{unipotent section} for $\SO (2, 1)$ one actually has that any orbit is either closed or dense. But this is false for $A$-orbits. A trivial example is the orbit $A.[e]$ of the identity coset which is a divergent orbit. Slightly less trivial is the example of an $A$-orbit of a point \[x \in \left [\begin{pmatrix} * & 0& 0 \\ 0& *& * \\ 0& *& * \end{pmatrix} \right]\] where essentially the action of $A$ degenerates to a rank-one action (one direction in $A$ acts in a trivial way sending every point to the cusp). The following example due to Shapira \cite{Shapira-Cassels-problem} of elements in $\SL (3, \R) / \SL (3, \Z)$ shows even this is not the only obstacle to $A.x$ being dense: Consider for any $\alpha \in \R$ the point
\begin{equation*}
p_\alpha= \left [\begin{pmatrix} 1 & 0& \alpha \\ 0& 1& \alpha \\ 0& 0& 1 \end{pmatrix} \right] 
.\end{equation*}
The $A$-orbit of $p_\alpha$ is certainly not $A$-periodic, but
\begin{equation*}
\overline {A . p _\alpha} \subset A . p _\alpha \cup \left [\begin{pmatrix} * & 0& * \\ 0& *& 0 \\ *& 0& * \end{pmatrix} \right] \cup \left [\begin{pmatrix} * & 0& 0 \\ 0& *& * \\ 0& *& * \end{pmatrix} \right];
\end{equation*}
see \cite{Lindenstrauss-Shapira} for details. Note the analogy to the possible behaviour allowed in Conjecture~\ref{torus conjecture}. Related examples of orbits of higher rank diagonal groups exhibiting this phenomena where given earlier by Macourant \cite{Maucourant-counterexample}, though not for a maximal diagonal group; another very interesting class of examples is investigated by Tomanov in \cite{Tomanov-counterexample}. 

For $x = [g] \in \SL (d, \R) / \SL (d, \Z)$, set
\begin{equation*}
\alpha _ 1 (x) =  \left (\inf_ {\mathbf n \in \Z ^ d \setminus 0} \norm {g \mathbf n} \right) ^{-1} \hspace{-10pt}
.\end{equation*}
The following conjecture seems to us plausible:

\begin{Conjecture}
Let $d \geq 3$, and let $x \in \SL (d, \R) / \SL (d, \Z)$ be such that
\begin{equation}\label{slow script rate}
\limsup_{a \in A} \frac {\log \alpha _ 1 (a.x) }{ \log \norm a} = 0
.\end{equation}
Then $\overline {A.x}=L.x$ for $A \leq L \leq \SL (d, \Z)$ and moreover $L.x$ is a periodic $L$-orbit (i.e. has finite volume).
\end{Conjecture}
As explained to us by Breuillard and Nicolas de Saxce \cite{Breuillard-de-Saxce-private}, the Strong Subspace Theorem of Schmidt \cite[\S VI.3]{Wolfgang-Schmidt-Diophantine-approximation} implies that if $g \in \SL(d, \overline {\Q})$ then \eqref{slow script rate} holds for~$[g]$ unless $g$ is in a proper $\Q$-parabolic subgroup of $\SL (d, \R)$. In particular we conjecture that if $g \in \SL(d, \overline {\Q})$, not contained in any proper $\Q$-parabolic subgroup of $\SL (d, \R)$, then $\overline {A.[g]}$ is homogeneous.

\medskip

Despite these difficulties, there are some positive results (not only conjectures) about orbit closures in this case. The first result in this direction is arguably Cassels and Swinnerton-Dyer result from their farsighted paper \cite{Cassels-Swinnerton-Dyer} that we already mentioned. In this paper,  Cassels and Swinnerton-Dyer prove that for the full diagonal group in $A<\SL (d, \R)$, every $A$-orbit $A . x$ that is itself non-periodic, but so that its closure $\overline{A.X}$ contains a periodic $A$-orbit, is unbounded. This allowed them to prove that Littlewood's conjecture (Conjecture~\ref{Littlewood Conjecture}) follows from Conjecture~\ref{C-SD conjecture}.

Using Ratner's Orbit Closure Theorem, Barak Weiss and the author were able to strengthen this as follows:

\begin{Theorem} [Weiss and L. \cite{Lindenstrauss-Barak}] \label{theorem with Barak}
Let $A . x$ be an orbits of the full diagonal group $A$ in $\SL (d, \R) / \SL (d, \Z)$ suppose that $\overline {A . x} \supset A . x _ 0$ with $A . x _ 0$ periodic. Then $\overline {A . x}$ is a periodic orbit of some group $L$ with $A \leq L \leq \SL (d, \R)$.
\end{Theorem}

\smallskip

An analogous result to Theorem~\ref{theorem with Barak} for $\SL (2, \Q _ p) \times \SL (2, \Q _ q) / \Gamma$ for $\Gamma$ an irreducible lattice arising from a quaternion division algebra was established earlier by Mozes \cite{Mozes-quaternions}; Mozes work is completely self-contained.
Theorem~\ref{theorem with Barak} was extended to inner forms of $\SL(d)$ (i.e. lattices arising from central simple algebras over $\Q$) by Tomanov in \cite{Tomanov-maximal-tori}.

We already noted that deciding e.g. if for $\alpha=\sqrt[3]{2}$
\begin{equation*}
\left\{ \begin{pmatrix} e ^ t& \\& e ^ {- t} \end{pmatrix} \left [\begin{pmatrix}  1& \alpha \\0& 1 \end{pmatrix} \right] : t \geq 0 \right\}
\end{equation*}
is bounded in $G / \Gamma = \SL (2, \R) / \SL (2, \Z)$ is a notoriously difficult question. Indeed, despite the fact that it is conjectured that for any irrational algebraic number $\alpha$ of degree $\geq 3$ the orbit above should be unbounded, not a single example of such an $\alpha$ is known.
For higher rank, e.g. for $\SL (3, \R) / \SL (3, \Z)$, one can at least give examples of explicit $A$-orbits of algebraic points that are dense. Indeed, Shapira and the author show in \cite{Shapira-Cassels-problem} and \cite{Lindenstrauss-Shapira} that if $\alpha, \beta$ are such that $1, \alpha, \beta$ span over $\Q$ a number field of degree 3 over $\Q$ then
\begin{equation*}
\overline{
A.\left [\begin{pmatrix} 1 & 0& \alpha \\ 0& 1& \beta \\ 0& 0& 1 \end{pmatrix} \right]^{\vphantom{l}}\,
}
 = \SL (3, \R) / \SL (3, \Z)
,\end{equation*}
This is related to another old result of Cassels and Swinnerton-Dyer, that showed in \cite{Cassels-Swinnerton-Dyer} that Littlewood's Conjecture (Conjecture~\ref{Littlewood Conjecture}) holds for such $\alpha,\beta$.

The strongest result to-date regarding Conjecture~\ref{C-SD conjecture} for general points is due to Einsiedler, Katok, and the author \cite{Einsiedler-Katok-Lindenstrauss} where using the classification of $A$-invariant measures of positive entropy on $\SL (d, \R) / \SL (d, \Z)$ which we obtained in that paper it was shown that for $d\geq 3$
\begin{equation*}
\dim_H \left\{ x \in \SL (d, \R) / \SL (d, \Z): \text{$A.x$ is bounded} \right\} = d-1
\end{equation*}
which implies that transverse to the flow direction (i.e. to $A$) the set of $x$ with a bounded $A$-orbit has zero Hausdorff dimension. Conjecturally, of course, this set is supposed to be a countable union of periodic $A$-orbits.

Finally we mention that Tomanov and Weiss \cite{Tomanov-Weiss} classified all closed $A$-orbits in $\SL (n, \R) / \SL (n, \Z)$ and more generally maximally split tori in arithmetic quotients, showing that an $A$-orbit $A.[g]$ is closed if and only if there is a $\R$-split maximal $\Q$-torus $\mathbb T < \SL (n)$ so that $A.[g]=g [\mathbb T (\R)]$. Their work builds upon the result of Margulis classifying all divergent $A$-orbits in $\SL (n, \R) / \SL (n, \Z)$ (such orbits correspond to $\Q$-split maximal $\Q$-tori). This has the striking consequence that if $F$ is a product of $n$ linearly independent linear forms in $n$-variables then $F(\Z ^ n)$ is discrete iff $F$ is proportional to an integral form \cite{Tomanov-values}.

\medskip

We now turn our attention to the question whether an analogue to the theorem of Mozes and Shah (Theorem~\ref{Mozes-Shah theorem}) holds for higher rank diagonal groups, where it seems the answer is mostly negative (but see \S\ref{sec: applications 1} for some significant positive results!)

For example, there are explicit examples of sequences of $A$-periodic orbits $A . x _ i$ in $\SL (d, \R) / \SL (d, \Z)$ for $A $ the $(d-1)$-dimensional diagonal group in $\SL (d, \R)$ so that the corresponding measures $m _ {A . x _ i}$ on $\SL (d, \R) / \SL (d, \Z)$ have escape of mass: there is a $0 \leq c<1$ so that for any compact $K$ for all $\epsilon > 0$ and all large enough $i$ we have that $\mu _ {A . x _ i} (K) < c + \epsilon$, and it is even possible to give such examples with $c=0$. Examples of $A$-periodic trajectories with escape of mass were noted in \cite{Einsiedler-Lindenstrauss-Michel-Venkatesh} (following a suggestion by Sarnak), with more elaborate examples (in particular with $c=0$) given by Shapira \cite{Shapira-full-escape} and David and Shapira \cite{David-Shapira}; implicitly these examples feature already in old work of Cassels \cite{Cassels-product}.
Escape of mass can also occur for a sequence of periodic measures for unipotent groups, or more generally a sequence of periodic measures that can arise as ergodic measures for unipotent groups, \emph{but only if the support of these measures, in its entirety, escape to infinity}, i.e. if we denote the sequence of measures by $\mu_i$ then for every compact set $K \subset G / \Gamma$ for every $i$ large enough, $K \cap \supp \mu _ i = \emptyset$. For the periodic $A$-orbits considered above there is a fixed set $K \subset \SL (d, \R) /  \SL (d, \Z)$ intersecting every one of them, indeed intersecting every $A$-orbit whether periodic or not.

Furthermore, assuming the equidistribution results of \cite{Einsiedler-Lindenstrauss-Michel-Venkatesh-III} hold in a quantitative way with polynomial error rates (which they surely should!), one can construct sequences of $A$-periodic orbits $A.x_i$ in $\SL (3, \R) / \SL (3, \Z)$ with volumes $\to \infty$ which converge weak$^*$ to a probability measure that gives positive mass to periodic orbit $A.y$ distinct from all the $A.x_i$.

\medskip

We end this section with a conjecture analogous to Conjecture~\ref{torus conjecture}.

\begin{Conjecture}\label{general orbit conjecture} Let $\mathbb G$ be an algebraic group over $\Q$, $S$ a finite set of places for $\Q$ containing $\infty$, and $G/\Gamma$ a corresponding $S$-arithmetic quotient saturated by unipotents (cf.\ Definition~\ref{saturated by unipotents}).
Let $A<G$ be a closed subgroup consisting of elements of class-$\mathcal A$ so that the projection of $A$ to $\prod_ {v \in S} (\mathbb G / \mathbb H) (\Q _ v)$ for any proper normal $\Q$-subgroup  $\mathbb H \lhd \mathbb G$ is of higher rank. Then for any $x \in G/\Gamma$ either $A.x$ is dense in $G/\Gamma$ or there are finitely many proper $\Q$-subgroups $\mathbb L_i < \mathbb G$ and $g_i \in G/\Gamma$ such that
\begin{equation*}
\overline {A.x} \subset A.x \cup \bigcup_ i  g_i[\mathbb L_i (\R)]
.\end{equation*}
\end{Conjecture}

\section{Applications regarding integer points and $\Q$-tori}\label{sec: applications 1}
The study of integer points on varieties is arguably the most basic problem in number theory. It seems at first sight rather surprising that the rigidity results for diagonalizable groups listed above could be relevant for such a problem. Fortunately they are, and perhaps a good point to start the discussion of this topic is by going back to the remarkable work of Linnik on the distribution of integer solutions to ternary quadratic equations, work which is presented in his book with the apt title ``Ergodic properties of number fields'' \cite{Linnik-book}, but in fact Linnik's farsighted work in this direction started even earlier in the late 1930s.

\subsection{Linnik's ergodic method for studying ternary quadratic forms using a one-parameter diagonalizable action}\label{Linnik's method}
Linnik considered several related problems: local to global results regarding which integers can be represented by an integer quadratic form in three variables, the distribution of integer points on a two dimensional sphere of radius $\sqrt m$ for $m \not\equiv 0,4,7 \bmod 8$, as well as the analogous problem regarding distribution of integer points on one and two sheeted hyperboloids in 3-space.

Consider in particular the distribution of integer points on the hyperboloid
\begin{equation*}
V _ d = \left\{ (a, b, c): b^2-4ac=d \right\}
\end{equation*}
where $d$ is an integer. Let $V _ d (\Z)$ denote the integer points on $V _ d$; these correspond to integral quadratic forms $a x ^2 + bxy + c y ^2$ of discriminant $d  = b ^2 -4 a c$, and let $V^*_d(\Z) \subseteq V _ d (\Z)$ the set of primitive points (i.e. triplets $(a,b,c)$ with no nontrivial common denominator). Note that for $V^*_d(\Z)$ to be nonempty, $d$ has to be $\equiv 0,1 \pmod 4$. The discriminant $d$ is said to be a fundamental discriminant if $V^*_d(\Z) = V _ d (\Z)$, i.e.\ if either $d$ is square-free and $\equiv 1 \pmod 4$ or $d = 4m$ with $m$ square-free satisfying $m\equiv 2 \text{ or } 3\pmod 4$. The action of $\GL (2)$ on binary quadratic forms gives us a natural action of $\GL (2, \Z)$ on $V ^*_ d (\Z)$ for every integer $d$.
 It is classical that $V^* _ d (\Z)$ consists of finitely many $\GL (2, \Z)$-orbits; indeed it is one of Gauss' remarkable discoveries that for a given $d$ one can define a natural commutative group law (which in this survey we denote by $\odot$) on $V^* _ d (\Z)/\GL (2, \Z)$; one way to characterise this group law is that if $[q_1] \odot [q_2] =[q_3]$ with $q_i \in V^* _ d (\Z)$ there are bilinear integral forms $\alpha,\beta$ so that 
 \begin{equation}\label{Gauss equation}
 q_1(n,m) q_2(l,s) = q_3(\alpha(n,m;l,s),\beta(n,m;l,s));
 \end{equation}
 cf.~\cite[\S1.3]{Cox-book}.
 As the identity in this group we take the $\GL (2, \Z)$-coset $[q_e]$ of $q_e=x^2-d' y^2$ if $d=4d'$ or $q_e=x^2+xy-(d-1)/4y^2$ if $d\equiv 1 \pmod 4$. For example, if $d=4d'$ equation \eqref{Gauss equation} applied to the triple $[q_e] \odot [q_e]=[q_e]$ is given explicitly by
 \[
 (n^2-d' m^2)(l^2-d' s^2) = (nl +d' ms)^2-d'(ns+ml)^2
 .\]

 \medskip
 
The following natural problem is a special case of an important class of counting questions raised by Linnik:
\begin{Question}\label{Linnik's question} Let $d _ i \to +\infty$. Let $\tilde V _ {d_i}$ be the sets $d _ i ^ {-1/2} V^* _ {d_i} (\Z) \subset V_1$. How are the points in these sets distributed? Let $m_{1}$ denote the unique (up to scalar) $\SL (2, \R)$-invariant measure on $V _ 1$. Do the points in $\tilde V _ {d_i}$ equidistribute, in the sense that for any nice subsets $E_1,E_2 \subset V_1$ (e.g. bounded open sets with $m _ 1$-null boundary)
\begin{equation*}
\frac {\#(\tilde V _ {d _ i} \cap E _ 1) }{ \#(\tilde V _ {d _ i} \cap E _ 2)} \to \frac {m_1 (E _ 1) }{ m _ 1 (E _ 2)}
?\end{equation*}
Similarly, let $d _ i \to {-\infty}$, and let $\tilde V _ {d _ i} = \absolute {d _ i} ^ {-1/2} V _ {d _ i} ^*(\Z) \subset V _ {-1}$. Do the points in $\tilde V_ {d _ i}$ become equidistributed in $V _ {-1}$?
\end{Question}

The answer to the question is YES, for any $d_i\to\infty$. This was proved by Duke \cite{Duke-hyperbolic-distribution} (building upon work of Iwaniec \cite{Iwaniec-half-integral-weight}), at least when $d_i$ is a sequence of fundamental discriminants  (which is, as implied by the name, the most fundamental [and hardest] case). Duke's proof is quantitative, and relies on estimates of Fourier coefficients of half integral weight Maass forms. Under an additional congruence condition on the sequence $d_i$, namely that there is some fixed prime $p$ so that $d_i$ are quadratic residues mod $p$ for all $i$ (i.e.\ $\bigl (\frac{ d_i }{ p }\bigr) = 1$), this equidistribution result was proved much earlier by Linnik and Skubenko \cite{Linnik-book, Skubenko-one-sheet}. In fact if $d_i \to +\infty$ a variant of Linnik's argument can be used to establish equidistribution with no (additional) side condition, as was shown by Einsiedler, Michel, Venkatesh and the author in \cite{Einsiedler-Lindenstrauss-Michel-Venkatesh-II}.

Set $G=\PGL (2, \R)$, $A<G$ the group of diagonal matrices and $\Gamma=\PGL (2, \Z)$. Note that $A$ is the stabilizer\footnote{Technically we are being slightly imprecise here, as $\PGL(2,\R)$ acts only on proportionality classes of quadratic forms.} of the quadratic form $q= x y$, and using this we can view $V_1$ as $A\backslash \PGL (2, \R)$. 
The above equidistribution question regarding $\tilde V _ {d_i}$, \ $d_i>0$, which we recall is a finite union of $\Gamma$-orbits in $V_1 \cong A \backslash G$, can be recast as a question regarding equidistribution of finite collections of closed $A$-orbits in $G/\Gamma$ as follows.
Consider a quadratic form $q (x,y) = ax^2 + b xy + cy^2 = a (x-\xi _ 1 y)(x-\xi _ 2 y)$ with $\xi _ {1,2} = \frac {-b \pm \sqrt d }{ 2 a}$. 
 We associates to $q$ the $A$-orbit $A.p_{a,b,c}$ in $G$ where
\begin{equation*}
p_{a,b,c}:=\begin{pmatrix} a&\frac{-b + \sqrt d }{ 2} \\ a&\frac{-b - \sqrt d} {2} \end{pmatrix}
,\end{equation*}
and where as before $A$ is the (one parameter) diagonal subgroup of $ G$. Clearly if $\gamma \in \Gamma$ the $A$-orbit corresponding to $q \circ \gamma$ will be $A.p_{a,b,c}\gamma$ and vice versa: if $A.p_{a,b,c}=A.p_{a',b',c'} \gamma$ for $(a,b,c),(a',b',c') \in V _ d^*(\Z)$ and $\gamma \in\Gamma$ then the corresponding quadratic forms $q$, $q '$ satisfy that $q = q ' \circ \gamma$. 

It follows that to each $\PGL (2, \Z)$-orbit in $V _ d  (\Z)$ there corresponds an $A$-orbit in $G/\Gamma$. A standard duality argument can be used to show (at least for $d_i$ that are not perfect squares, though the case of $d_i$ perfect squares can also be handled this way) that equidistribution of the sequence of sets $\tilde V _ {d_i}$ in $V_1 \cong A \backslash G$ in the sense of Question~\ref{Linnik's question} is equivalent to the equidistribution of the sequence of collections of closed $A$-orbits
\[
\mathcal T _ {d_i} = \left\{ A .p_{a,b,c}: (a, b, c) \in V _ {d _ i} ^{*} (\Z) \right\}
\]
in $G/\Gamma$, i.e. that for every $f,g \in C _ 0 (G/\Gamma)$ with $\int g\,dm\neq 0$
\begin{equation}\label{Linnik equidistribution}
\frac {\displaystyle \sum_ {A.p \,\in \mathcal T _ {d_i}}\textstyle \int_ {A.p} f }{ \displaystyle \sum_ {A.p \,\in \mathcal T _ {d_i}} \textstyle\int_ {A.p} g} \to \frac{\int_ {G/\Gamma} f \,d m}{\int_ {G/\Gamma} g \,d m}\hspace{0.18cm}.
\end{equation}
Note that we present the equidistribution in the above form to allow for sequences $d_i$ that contain perfect squares, as in that case the $A$-orbits $A .p_{a,b,c}$ are divergent. For $d_i$ a sequence avoiding perfect squares one can take $g\equiv 1$ instead.

\label{two more ways page}
We give two more ways to look at the points $ V^*_d (\Z)/\GL(2,\Z)$ which will be important for us later.  
\begin{enumerate}[label=\textbf{\Roman*.},ref=\textbf{\Roman*}]
\item \label{isomorphism to class group}
The cosets $\left\{ [p _ {a, b, c}] _ \Gamma : (a,b,c) \in V _ d ^{*} (\Z) \right\}$ correspond in an obvious way to the $\Z$-modules $\bar p _ {a, b, c}$ spanned by $a$ and $\frac {- b + \sqrt d }{ 2}$ which (for $d$ not a perfect square) are in fact ideals in the order $\mathcal{O} _ d$ of discriminant $d$ (for fundamental discriminants, $\mathcal{O} _ d$ is the ring of integers in $\Q (\sqrt d)$; if $d=f^2d'$ with $d'$ fundamental $\mathcal{O} _ d$ is a subring of $\mathcal{O} _ {d '}$ containing the identity of index $f$ in $\mathcal{O} _ {d '}$). For $(a,b,c),(a',b',c') \in V _ d^*(\Z)$ we have that $[p _ {a', b', c'}] _ \Gamma \in A.[p _ {a, b, c}] _ \Gamma$ if and only if the ideals $\overline p _ {a', b', c'}$ and $\overline p _ {a, b, c}$ are in the same ideal class, i.e. for some $k \in \Q (\sqrt d)$ we have that $(k \mathcal{O} _ d) \cdot \overline p _ {a', b', c'} = \overline p _ {a, b, c}$. An observation that can be attributed essentially to Dirichlet is that the Gauss composition law on $\GL(2,\Z)$-cosets in $ V^*_d (\Z)$ is the same group law as the group law in the ideal class group $\cl (\mathcal O_d)$ of the order $\mathcal O_d$.

\item \label{isomorphism to adelic torus}
For any positive $d \in \Z$, and $(a,b,c) \in V^*_d (\Z)$, the group $p _ {a,b,c} ^{-1} A p _ {a,b,c}$ is the group of $\R$-points of a $\Q$-torus $\mathbb T_{a,b,c} < \PGL (2)$. Moreover, it is not hard to see that $\mathbb T_{a,b,c}$ is $\Q$-split iff $d$ is a perfect square. It follows that the orbit $A.[p_{a,b,c}]$ is a closed $A$-orbit in $G/\Gamma$, and this $A$-orbit $A.[p_{a,b,c}]$ is a periodic $A$-orbit if $d$ is not a perfect square and divergent otherwise.
The tori $\mathbb T_{a,b,c}$ are conjugate to each other in $\PGL(2)$ over $ \Q$ but not over $\Z$:  for $(a,b,c),(a',b',c') \in V^*_d(\Z)$ the tori $\mathbb T_{a,b,c}$ and $\mathbb T_{a',b',c'}$ are conjugate to each other over $\Z$ if and only if $p_{a',b',c'} \in p_{a,b,c}  \Gamma$.

These collections $\mathcal T _ {d }$ of $A$-orbits in $G / \Gamma$ can be described very succinctly in the language of the adeles.
Let $\mathbb A$ denotes the adele ring of $\Q$, and let \[\pi: \PGL (2, \mathbb A) / \PGL (2, \Q) \to \PGL (2, \R) / \PGL (2, \Z)=G/\Gamma\] be the natural projection, which takes the coset $[(g,g_2,g_3, \dots)]_{\PGL (2, \Q) }$ of an element $(g,g_2,g_3, \dots) \in \PGL (2, \mathbb A) $ with $g \in \PGL (2, \R)$ and $g_p \in \PGL (2, \Z _ p)$ for every (finite) prime $p$ to the coset $[g]_\Gamma$ in $G / \Gamma$.
Then for any $(a, b, c) \in V ^{*} _ d (\Z)$
\begin{equation*}
\mathcal T _ {d} = p_{a,b,c} \, \pi\! \left ( [\T_{a,b,c}(\mathbb A)] \right);
\end{equation*}
since this is valid for any choice of $(a,b,c)$ we may as well take the explicit choice of $(a,b,c)=(1,0,-d/4)$ (for $4|d$) or $(a,b,c)=(1,1,-(d-1)/4)$ (for $d\equiv 1 \pmod 4$).
\end{enumerate}
Up to minor changes --- replacing $A$ with the compact group $K = \left\{ \begin{pmatrix} \cos \theta& \sin \theta \\ - \sin \theta& \cos \theta \end{pmatrix} \right\}$ and taking $p_{a,b,c}$ to be $\begin{pmatrix} a&-b/2 \\ 0& \sqrt {\absolute{d}}/2 \end{pmatrix}$ --- the correspondences described in \ref{isomorphism to class group} and \ref{isomorphism to adelic torus} above hold also for $d$ negative. Of course, $A$ and $K$ are quite different $\R$-groups, with $K$ being $\R$-anisotropic and compact while $A$ is $\R$-split.

From the remainder of this subsection we restrict our attention to $d_i>0$ not perfect squares (i.e. $\T_{a,b,c}$ \ $\Q$-anisotropic); for a discussion of the isotropic case see \cite{Oh-Shah-divergent, David-Shapira-divergent, Shapira-Zheng-divergent}.
As explained in \cite{Einsiedler-Lindenstrauss-Michel-Venkatesh-II}, in modern terminology Linnik's approach can be interpreted as the following 3-step strategy:
\begin{enumerate}[label=\textbf{\Alph*.},ref=\textbf{\Alph*}]
\item \label{first step for Linnik} Let $\mu _ i$ be the probability measure given by
\begin{equation*}\mu_i(f) =
\frac {\displaystyle \sum_ {A.p \,\in \mathcal T _ {d_i}} \textstyle \int_ {A.p} f }{ \displaystyle \sum_ {A.p \,\in \mathcal T _ {d_i}} \textstyle \int_ {A.p} 1}
\end{equation*}
One needs to establish that this sequence of measures is \emph{tight}, i.e. that for every $\delta>0$ there is a compact $X _\delta \subset G / \Gamma$ so that $\mu _ i (X _\delta) > 1 - \delta$ for all $i$. Linnik establishes this via analytic number theory, in a way that is closely related to a key step in \cite{Einsiedler-Lindenstrauss-Michel-Venkatesh-III} which we discuss below, but this can also be established using purely ergodic theoretic means (cf. \cite{Einsiedler-Lindenstrauss-Michel-Venkatesh-II}).

\item 
Then Linnik proves an upper bound on the measure of small tubes transverse to the $A$ action of radius $r\geq d_i^{1/4}$, on average:
\begin{equation}\label{Linnik's estimate}
\int_ {X_\delta} \mu _ i (B(r,1,x)) \,d \mu _ i (x) \ll_{\delta, \epsilon} r^{2} d_i ^ \epsilon
\end{equation}
for $\epsilon, \delta > 0$, $X_\delta$ as above, and
\[
B(r,1,x') = \left\{ \begin{pmatrix}  1& s \\0& 1 \end{pmatrix} \begin{pmatrix} e ^ t& \\& e ^ {- t} \end{pmatrix} \begin{pmatrix}  1& 0 \\s'& 1 \end{pmatrix}. x': \absolute s, \absolute {s'} < r, \absolute t <1 \right\}.
\]
 An important point here is that the $d_i^\epsilon$ term in the right hand side implies that this estimate is meaningful only for $r>d_i^{-c}$, hence for each $i$ we obtain information regarding the distribution of $\mu _ i$ at a different scale.

The estimate \eqref{Linnik's estimate}, which Linnik called the Basic Lemma, is key to the whole approach. note that the exponent $2$ in the right inside of \eqref{Linnik's estimate} is \emph{sharp}. It is a deep bound, that relies on results of Siegel and Venkov on quadratic forms, and is closely related to the Siegel Mass Formula (cf.\ \cite[Appendix A]{Einsiedler-Lindenstrauss-Michel-Venkatesh-II} for a self-contained treatment).

\item\label{third step for Linnik} Now somehow one needs to upgrade the sharp non-concentration (on average) estimate \eqref{Linnik's estimate} to an equidistribution statement: to both lower and upper bounds on the measure of fixed sized subsets of $G / \Gamma$. One way to proceed, explained in \cite{Einsiedler-Lindenstrauss-Michel-Venkatesh-II} is as follows: by passing to a subsequence if necessary, and in view of the tightness of the sequence of measures $\mu _ i$ discussed in step~\ref{first step for Linnik}, we can assume that $\mu _ i$ converges to some $A$-invariant measure $\mu$ as $i \to \infty$ in the weak$^*$-topology.
The action of $A$ on $G/\Gamma$ is a prime example of a rank-one diagonalizable group actions, one which does \textbf{not} satisfy the type of rigidity provided by Theorem~\ref{split theorem semisimple} of the other measure classification theorems discussed in \S\ref{higher rank measure rigidity}. Because there are so many $A$-invariant measures in $G / \Gamma$ it is in general a very hard to prove that a limiting measure obtained from a number theoretic construction will be the uniform measure (cf.\ the discussion at the beginning of \S\ref{torus orbits section}). However the fact that the estimate \eqref{Linnik's estimate} is \emph{sharp} rescues us, as \eqref{Linnik's estimate} together with the subadditivity of entropy allows us to deduce that the entropy of $\mu$ with respect to the action of $a_t \in A$ is maximal, and on $G / \Gamma$ there is a \emph{unique measure of maximal entropy}.

\end{enumerate}

Linnik and Skubenko did not quite follow the method outlined in step \ref{third step for Linnik}. To begin with, they consider the dynamics not for the diagonal subgroup in $\PGL  (2, \R)$ but to a diagonal group over $\Q_p$: Linnik and Skubenko assumed that for some fixed prime  $p$, the sequence $d_i$ satisfied $\bigl (\frac { d _ i}{p } \bigr) = 1$, which implies that the measures $\mu _ i$ (as well as any limiting measure $\mu)$ can be lifted to probability measures on $\PGL (2, \R) \times \PGL (2, \Q _ p) / \PGL (\Z [1/p])$ that are invariant under the diagonal subgroup of $\PGL (2, \Q _ p)$. This $p$-adic dynamics is symbolic in nature, which facilitated the analysis. Moreover they did not first pass to the limit, which allowed them to give rates of equidistribution (even if logarithmic rather than polynomial as in the work of Duke). 
A third alternative\footnote{Some may say this is more of a development of Linnik's original method.} to step \ref{third step for Linnik} using property $\tau$ (i.e.\ using spectral gaps) was given by Ellenberg, Michel and Venkatesh in \cite{Ellenberg-Michel-Venkatesh}, which gives another readable modern interpretation of the Linnik method, and in addition raises an important joining question to which they are able to give a partial answer (see below).

Linnik's method is not limited to the discriminant form $b^2-4ac$, but is applicable to any integral ternary quadratic form, in particular to the form $a^2+b^2+c^2$, i.e. to the distribution of integer points on the sphere. Both $b^2-4ac$ and $a^2+b^2+c^2$ are quadratic forms with the nice property that any other integral form that is equivalent to them over $\R$ and $\Z _ p$ for all $p$ is in fact equivalent to them over $\Z$. The collection of all integral quadratic forms equivalent  over $\R$ and $\Z _ p$ for all $p$ to a given quadratic form is called the \emph{genus} of the quadratic form\footnote{A somewhat anachronistic terminology as this genus has nothing to do with the genus of any surface.}.\label{genus page} For general integral quadratic forms one needs to study all forms in the genus in order to prove equidistribution of integer points on each of the corresponding quadratic surfaces. The form $a^2+b^2+c^2$ is treated explicitly by Linnik in his book \cite{Linnik-book} and earlier works, and is also the case explained in \cite{Ellenberg-Michel-Venkatesh}; the case of general quadratic forms is discussed e.g.\ in Linnik's paper \cite{Linnik-genera}. See also \cite{Wieser-Linnik} for a nice modern exposition by Wieser.

\subsection{Going beyond Linnik --- joint equidistribution using rigidity of joinings for higher-rank diagonalizable actions}\label{beyond Linnik section}

In the previous section we considered how points on the one or two sheeted hyperbolic
\begin{equation*}
V _ d (\Z) = \left\{ (a, b, c): b^2-4ac=d \right\}
\end{equation*}
project onto the unit one or two sheeted hyperboloid $V_1$ or $V_{-1}$  respectively depending on the sign of $d$, and similarly regarding projection of points on the sphere of radius $\sqrt {-d}$ (for notational convenience, we will use negative integers to parameterize spheres)
\begin{equation*}
S _ d (\Z) = \left\{ (a, b, c: \in \Z ^ 3: a ^2 + b ^2 + c ^2 = - d \right \}
\end{equation*}
projects to the unit sphere $\mathbb S$. As explained in \ref{isomorphism to class group} and \ref{isomorphism to adelic torus} on p.~\pageref{two more ways page} for the special case of the one sheeted hyperboloid (i.e. $V_d(\Z)$ for $d>0$; cf.\ also the paragraph immediately afterwards regarding the modification for $d<0$), these distribution problems regarding the integer points $V _ d (\Z)$ and $S_d(\Z)$ can be interpreted in terms of the ideal class group of $\Q(\sqrt{d})$ or translated into questions regarding the distribution of suitable translates of the adelic points of $\Q$-tori $\mathbb T_{d} (\mathbb A)$ in $\mathbb G (\mathbb A)/\mathbb G(\Q)$ for $\mathbb G$ being the $\Q$-group $\PGL (2)$ or $\SO (3)$ in the hyperboloid and sphere cases respectively.

We will be mainly interested in the harder case of $\mathbb T_{d}$ \ $\Q$-anisotropic, i.e. $d$ not a perfect square, in which case $\mathbb T_{d}$ will split over the quadratic extension $\Q(\sqrt{d})$ of~$\Q$. In particular, the $\Q$-torus $\mathbb T _ d$ will be split over $\R$ for $d > 0$, i.e.\ for one sheeted hyperboloids, and over $\Q_p$ iff $\bigl (\frac {d }{ p}\bigr)=1$. If one wants to follow a scheme as in \ref{first step for Linnik}-- \ref{third step for Linnik} to prove equidistribution using dynamical ideas, particularly if one follows~\ref{third step for Linnik} to construct a limiting measure $\mu$ out of a sequence $d_i \to \pm \infty$ and use dynamics to study this limiting measure, one needs to assume either that $d_i >0$ for all $i$ or that there is a fixed prime $p$ so that $\bigl (\frac {d_i }{ p}\bigr)=1$ for all $i$.

We can strengthen our assumptions, and require \textbf{two} places $v,w \in \left\{ \infty, \text{primes} \right\}$ at which the tori  $\mathbb T _ {d_i}$ splits, i.e. that $\bigl (\frac {d_i }{v}\bigr) = \bigl (\frac {d_i }{w}\bigr)= 1$ for all $i$; to allow also the case of $v$ or $w=\infty$ , we define
$\bigl (\frac {d }{\infty}\bigr) = {\rm sign}(d)$. If we make this assumption we will obtain a limiting measure $\mu$ on $\mathbb G (\mathbb A) / \mathbb G (\Q)$ (or, if we prefer, on an $S$-arithmetic quotient of $\mathbb G$ for any $S$ containing $v,w, \infty$) that is invariant under a \emph{higher rank} diagonal group, on which we can try to apply the measure rigidity theorems presented in \S\ref{measure rigidity section}, and in particular the joining classification theorem, Theorem~\ref{higher rank joining}.

\medskip

We now describe an arithmetic consequence of the rigidity of higher rank diagonal groups obtained in this way by Aka, Einsiedler, and Shapira \cite{Aka-Einsiedler- Shapira}. Let $d$ be a negative integer. By the Three Squares Theorem of Legendre and Gauss, $S_{d}(\Z)$ is non-empty, iff $d\not\equiv 1 \bmod 8$ (recall that in our conventions, $d$ is negative!)
Consider for any integer vector $\mathbf n \in S_{d_i} (\Z)$ the lattice in the plane orthogonal to $\mathbf n$ (with respect to the standard inner product on $\R^3$)
\[
\Lambda _ { \mathbf n} = \left\{ \mathbf x \in \Z ^ 3: \mathbf x \cdot \mathbf n = 0 \right\}
.\]
Let $\mathbf v _ 1, \mathbf v _ 2$ be generators of $\Lambda _ {\mathbf n}$ (considered as an additive group); then $\mathbf v _ 1, \mathbf v _ 2$ give rise to a positive definite binary quadratic form $q(x,y)=\norm {x \mathbf v _ 1 + y \mathbf v _ 2}^2$. The integer quadratic form $q$ will have (negative) discriminant $4d$ (we leave this as an exercise to the reader\dots), and given $\mathbf n$ the form $q$ is well defined up to the action of $\GL (2, \Z)$ on $V_{4d}(\Z)$. To be slightly more precise, $\mathbf n$ gives an orientation on the plane ${\mathbf n}^\perp$, so if we chose $\mathbf v _ 1, \mathbf v _ 2$ to be a basis with positive orientation the form $q$ is well defined up to the action of $\SL (2, \Z)$.
 Thus we get a map $\alpha: S _ {d} (\Z) \to V _ {4d} (\Z) / \SL (2, \Z)$. This map is neither injective nor onto, but is close to being both --- up to a bounded integer factor, both the kernel and co-kernel of this map is of size $2^{r-1}$, where $r$ is the number of distinct primes dividing $d$ (this is already less obvious, but was understood already by Gauss). 
 Let $m_Y$ denote the $\SL(2,\R)$-invariant measure on $Y=V_{-1}/ \SL (2, \Z)$ and $m_S$ the uniform measure on the unit sphere $\mathbb S$ (both normalized to be probability measures). We remark that $Y$ can naturally be identified with $\H / \PSL(2, \Z)$, with $\H$ the hyperbolic plane. 
The map $\alpha$ is close enough to being 1-1 that e.g.\ if one assumes that $d_i \to {-\infty}$ with $d_i \not\equiv 1 \bmod 8$ squarefree and $\bigl (\frac {d_i }{p}\bigr)=1$ for some fixed prime $p$ then it follows e.g.\ using Linnik's methods that the projections of $\left\{ \alpha (\mathbf n): \mathbf n \in S _ {d _ i} (\Z)\right\}$ to $Y$  becomes equidistributed with respect to $m_Y$. Recall also that Linnik showed under these conditions that the collection of points $S _ {d _ i}(\Z)$ projected to the unit sphere $\mathbb S$, become equidistributed with respect to $m_S$. Using Theorem~\ref{higher rank joining}, Aka, Einsiedler and Shapira where able to upgrade these two statements to a joint equidistribution statement:

\begin{Theorem} [Aka, Einsiedler and Shapira \cite{Aka-Einsiedler-Shapira}]\label{Aka-Einsiedler-Shapira theorem} Let $p,q$ be two distinct finite primes, and $d _ i \to {-\infty}$ a sequence of square free negative integers $\not \equiv 1 \bmod 8$ so that $\bigl (\frac {d_i }{p}\bigr) = \bigl (\frac {d_i }{q}\bigr)= 1$ for all $i$. Then the projection of the sets
\begin{equation*}
\left\{ (\mathbf n, \alpha (\mathbf n)): \mathbf n \in S _ {d _ i} \right\} \subset S _ {d _ i} (\Z) \times (V _ {4 d _ i}(\Z) / \SL (2, \Z))
\end{equation*}
to $\mathbb S  \times Y$ becomes equidistributed with respect to the measure $m_S \times m_Y$ on this space as $i \to \infty$.
\end{Theorem}

\noindent
In other words, for any nice subsets $E \subset \mathbb S $ and $F \subset Y $
\begin{equation*}
\frac {\#\left\{ \mathbf n \in S _ {d _ i} (\Z):  \text{$\absolute {d _ i} ^ {-1/2}\mathbf n \in E$ and $   \absolute {4d _ i} ^ {-1/2}\alpha (\mathbf N) \in F$} \right\} }{ \#S _ {d _ i} (\Z)} \to m_{S} (E) \cdot m_{Y }(F)
.\end{equation*}

Unlike the individual equidistribution on $\mathbb S$ and $V _ {-1} / \SL (2, \Z)$ which can also be proved using analytic number theoretic tools (indeed, the analytic tools give significantly sharper results), there does not seem to be a plausible approach using currently available technology to proving this joint equidistribution statement using the techniques of analytic number theory or automorphic forms. For more information in this direction we refer the reader to \cite[Appendix]{Aka-Einsiedler-Shapira-arxiv} by Ruixiang Zhang\footnote{\cite{Aka-Einsiedler-Shapira-arxiv} is the arXiv version of \cite{Aka-Einsiedler-Shapira}.}.

Theorem~\ref{Aka-Einsiedler-Shapira theorem} turns out to be closely related to the following equidistribution result stated in terms of the class group:

\begin{Theorem}\label{class group joining theorem} Let $d _ i \to {-\infty}$ be a sequence of negative integers $ \equiv 0$ or $1 \bmod 4$ so that there are two primes $p,q$ for which $\bigl (\frac {d_i }{p}\bigr) = \bigl (\frac {d_i }{q}\bigr)= 1$ for all $i$. Identifying as before elements of $V _ {d_i} ^{*} (\Z)$ with primitive integral quadratic forms, we fix for every $i$ an arbitrary integral form $q_i \in V _ {d_i} ^{*} (\Z)$, and define a collection of points $\tilde V_{d_i}^{(2)} \subset (V _ {d _ i} (\Z) / \GL (2, \Z)) ^2$ by
\begin{equation}\label{q, q squared equation}
\tilde V_{d_i}^{(2),q_i} = \left\{ ([q], [q_i]\odot [q]\odot [q]): q \in V _ {d_i} ^{*} (\Z) \right\}
.\end{equation}
Then the projection of these collections to  $Y ^2$ become equidistributed with respect to $m_Y \times m_Y$ as $i \to \infty$.
\end{Theorem}

\Head{Note} We restrict ourselves to $d_i<0$ for purely aesthetic reasons, as in this case the relation of the equidistribution statement to integer points is cleanest. In fact, taking $d_i \to +\infty$ is even better, as then only one additional split place is needed, i.e. one need only assume the existence of one prime $p$ for which $\bigl (\frac {d_i }{p}\bigr) = 1$ for all $i$.

\begin{proof}[Sketch of proof]
For simplicity of notations, assume $4|d_i$ for all $i$ (the modification to $d_i \equiv 1 \pmod 4$ poses no additional difficulties). Let $q_e = x^2-\frac{d_i}4 y^2$, and let $\mathbb T _{i}$ be the corresponding adelic torus as in p.~\pageref{two more ways page}, \ref{isomorphism to adelic torus} (adapted for $d < 0$). Concretely, we can take $\mathbb T _{i}$ to be the stabilizer of the proportionality class of $q_e$ in $\PGL (2, \R)$. The $\Q$-torus $\mathbb T _ i$ is anisotropic over $\Q$ and even over $\R$, but by our assumption on $\bigl (\frac {d_i }{p}\bigr)$ and $\bigl (\frac {d_i }{q}\bigr)$ will be split over $\Q _ p$ and $\Q _ q$. Let $G= \PGL (2, \R)$, $\Gamma = \PGL (2, \Z)$, $K$ the maximal compact subgroup of $\PGL (2, \Z)$ as in p.~\pageref{two more ways page}, $S = \left\{ \infty, p, q \right\}$, $G _ S = \prod_ {v \in S} \PGL (\Q _ v)$, $\Gamma _ S = \PGL (\Z[1/pq])$, diagonally embedded in $G_S$.

Let $\tilde V_{d_i} = \absolute {d _ i} ^ {-1/2} V ^{*}_{d_i} (\Z) \subset V_{-1}$. Consider the natural projections $\pi _ S, \pi_Y$ as in the diagram below
\begin{equation*}
\mathbb G (\mathbb A) / \mathbb (\Q) \xrightarrow{\pi _ S} G _ S / \Gamma _ S \xrightarrow {\pi_Y} K \backslash G / \Gamma \cong Y
,\end{equation*}
and let $\pi=\pi_Y \circ \pi_S$. We denote by $\pi'$, $\pi'_Y$ the unnormalized form of $\pi$ and $\pi_Y$, i.e. the corresponding maps to $V_{d}/\GL (2, \Z)$. For suitable choice of $g _ i \in \PGL (2, \mathbb A)$ we have that
\[ \pi (g_i [\T _ i (\mathbb A)]) = \tilde V_{d_i} / \GL (2, \Z)
\]
and that $\pi _ S (g _ i \T _ i (\mathbb A ))$ is invariant under the diagonal group $A_1 < \PGL (2,\Q _ p)$ and $A_2 < \PGL (2,\Q _ q)$.
In particular, there is a $t _ i \in \T _ i (\mathbb{A})$ so that $\pi (g _ i [t _ i])$ is the point in $Y $ that corresponds to the $\GL (2, \Z)$-coset of $q _ i$.

Reconciling the two points of view on the set $V _ d ^{*} (\Z)$ given in~\ref{isomorphism to class group} and~\ref{isomorphism to adelic torus} on p.~\pageref{two more ways page}, one verifies that
\begin{equation*}
\tilde V_{d_i}^{(2),q_i} = \left\{ ([q], [q_i]\odot [q]\odot [q]): q \in V _ {d_i} ^{*} (\Z) \right\}
\end{equation*}
is equal to
\begin{equation*}
\left\{ (\pi'(g_i[t]), \pi ' (g_i[t_i \cdot t \cdot t])): t \in \T _ i(\mathbb{A}) \right\},
\end{equation*}
in particular $\tilde V_{d_i}^{(2),q_i}$ is the projection of the $\left\{ (a_1a_2,a_1^2a_2^2):a_1 \in A_1, a_2 \in A_2 \right\}$-invariant subset
\begin{equation*}
\left\{ (\pi _ S (g_i[t]), \pi _ S (g_i[t_i t^2])): t \in \T _ i(\mathbb{A}) \right\} \leq (G _ S / \Gamma _ S) ^2
\end{equation*}
to $(V _ {d} / \GL (2, \Z)) ^2$.

The projection of the first coordinate in $\tilde V _ {d _ i} ^ {(2),q_i}$ to $Y$ equidistributes by the work of Linnik \cite{Linnik-book} or by Duke \cite{Duke-hyperbolic-distribution} (since we already assumed two split places --- $\bigl (\frac {d_i }{p}\bigr) = \bigl (\frac {d_i }{q}\bigr)= 1$ for all $i$ --- we may as well use Linnik who needs only one).
The second coordinate in $\tilde V _ {d _ i} ^ {(2), q _ i}$ does not run over all of $V _ {d _ i} ^{*} (\Z) / \GL (2, \Z)$ but rather over a sub-collection, say $V^{\rm second}_{d_i}$ of index equal to the 2-torsion in the class group $\cl(\mathcal{O} _ {d _ i})$ of $\mathcal{O} _ {d _ i}$.

Fortunately, already Gauss understood the 2-torsion in the class group of quadratic fields (remarkably, even today we do not understand 2-torsion of the class group in fields of higher degree!), and its size is nicely controlled by the number of divisors of $d _ i$; in particular has size $\ll \absolute { d_i}^\epsilon$ for all $\epsilon > 0$ (by a theorem of Siegel the size of  $\cl(\mathcal{O} _ {d _ i})$ is (noneffectively) $\gg \absolute { d_i}^{1/2 - \epsilon}$ for all $\epsilon$).

To prove the equidistribution of $V_{d_i} ^ {\rm second}$, one can either quote a result of Harcos and Michel \cite{Harcos-Michel} that can be viewed as an extension of Duke's work, or use ergodic theory: $V_{d_i} ^ {\rm second}$ is the projection of a $\left\{ a_1^2a_2^2: a _ 1 \in A _ 1, a _ 2 \in A _ 2 \right\}$-invariant subset of $G _ S / \Gamma _ S$ that can be treated using Linnik's method as the analogue of \eqref{Linnik's estimate} will holds also for $V_{d_i} ^ {\rm second}$. For more details\footnote{At least for the case of fundamental discriminants, though the general case is not more complicated.}, we refer the reader to \cite[\S4]{Aka-Einsiedler-Shapira}.

Once the equidistribution of each component of $\tilde V ^ {(2), q _ i} _ {d _ i}$ separately has been proved, Theorem~\ref{higher rank joining} takes care of the rest. A key point is that there is no nontrivial algebraic joining in $(G _ S / \Gamma _ S) ^2$ invariant under $\left\{ (a_1a_2,a_1^2a_2^2):a_1 \in A_1, a_2 \in A_2 \right\}$.
\end{proof}

It is a folklore conjecture that for any integer $k$, the $k$-torsion in $\cl(\mathcal{O} _ {d})$ is $\ll \absolute d ^ \epsilon$ as $\absolute d \to \infty$ (see e.g.~\cite{Ellenberg-Venkatesh-torsion}). \textbf{Assuming this conjecture} the method of Theorem~\ref{class group joining theorem} would give that for any $k$, any sequence of negative integers $d _ i \to {-\infty}$ with  $d_i \equiv 0$ or $1 \bmod 4$ and $\bigl (\frac {d_i }{p}\bigr) = \bigl (\frac {d_i }{q}\bigr)= 1$ for two fixed primes $p$,$q$, and any choice of $q_{i,1},q_{i,2},\dots  q_{i,k} \in V _ {d_i} ^{*} (\Z)$ we have that the projection of
\[
\left\{ ([q_{i,1}]\odot [q], [q_{i,2}]\odot [q]^{\odot 2},\dots, [q_{i,k}]\odot [q]^{\odot k} ): q \in V _ {d_i} ^{*} (\Z) \right\}
\]
to $Y ^ k$ becomes equidistributed.

Consider now another collection of points in $(V _ {d_i} / \GL (2, \Z)) ^2$:
\begin{equation}\label{definition of join collection}
\tilde V_{d_i}^{{\rm join},q_i}=\left\{ ([q], [q_{i}]\odot [q]: q \in V _ {d_i} ^{*} (\Z) \right\}
\end{equation}
depending on the discriminant $d _ i$ as well as $q _ i \in V _ {d _ i} ^{*} (\Z)$.

It would seem to be a simpler collection to study than the ``non-linear'' collection $\tilde V_{d_i}^{(2),q_i}$ defined by \eqref{q, q squared equation}, if only because the equidistribution of each of the two projections to $Y$ follows literally from Duke's theorem. However, this intuition turns out to be misguided: studying the distribution of the collections $\tilde V_{d_i}^{{\rm join},q_i}$ turns out to be substantially subtler than that of the collection $\tilde V_{d_i}^{(2),q_i}$ for a simple reason: after passing to a limit, \textbf{there are nontrivial joinings} that need to be ruled out. And indeed, without further assumptions the projections of collection of points $\tilde V_{d_i}^{{\rm join},q_i}$ to $Y\times Y$ need not equidistribute.
For instance, if one takes $q_i$ to be $x^2-(d_i/4)y^2$ or $x^2+xy-\frac{(d_i-1)}4y^2$ (so that $[q_i]$ is the identity for Gauss' group law on $V _ {d_i} ^{*} (\Z) /\GL (2, \Z)$) the projection of this collection of points to $Y ^2$ is supported on the diagonal $ \left\{ ([v],[v]): v \in  Y \right\}$ hence certainly does not equidistribute!

A similar problem holds if $[q_i]$ has small ``size'' $\operatorname{\frak{N}} ([q_i])$. To define $\operatorname{\frak{N}} (\bullet)$ we use the correspondence in~\ref{isomorphism to class group} on p.~\pageref{two more ways page} between elements of $\tilde V ^{*} _ {d _ i} (\Z)/\GL (2, \Z)$  and elements of the ideal class group  $\cl(\mathcal{O} _ {d _ i})$. If $I \lhd \mathcal{O} _ {d _ i}$ is the ideal corresponding to the quadratic form $q \in \tilde V ^{*} _ {d _ i} (\Z)/\GL (2, \Z)$ we define
\begin{equation*}
\operatorname{\frak{N}} ([q]) = \min \left\{ \operatorname{\frak{N}} (J): J \lhd \mathcal O_{d_i}, J \sim I \right\}
.\end{equation*}
If we consider a sequence $d _ i \to {-\infty}$ and $q _ i \in V ^{*} _ {d_i} (\Z)$ with $\operatorname{\frak{N}} ([q _ i])$ bounded (say less than $N$), for similar reasons the collections $\tilde V _ {d _ i} ^ {{\rm join}, q _ i}$ will be restricted to a subset of $Y \times Y$ of dimension $\dim Y = 2$: the union of the graphs of the Hecke correspondences on $Y \times Y$ of order $\leq N$.

Michel and Venkatesh conjectured this is the only obstruction:
\begin{Conjecture} [Michel and Venkatesh {\cite [Conj.~2]{Michel-Venkatesh-ICM}}]\label{joining conjecture}
let $d _ i \to {-\infty}$ along the sequence of fundamental discriminants. For each $i$, let $q_i \in V ^{*} _ {d _ i} (\Z)$, and assume that $\operatorname{\frak{N}} ([q _ i]) \to \infty$. Then the projection of the collection $\tilde V_{d_i}^{{\rm join},q_i}$
to $Y \times Y$ equidistribute as $i \to \infty$.
\end{Conjecture}

We recall that for $d < 0$ (ineffectively) the size of $\cl (\mathcal{O} _ {d})$, which as we have seen can be identified with $V_ {d} ^{*} (\Z) / \GL (2, \Z)$,  is $\absolute {d}^{1/2+o(1)}$. The number of ideal classes $[I]$ in $\cl (\mathcal{O} _ {d})$ with $\operatorname{\frak{N}} ([I]) < N$ can be easily seen to be $\ll N^{1+o(1)}$. Moreover, by a simple application of Minkowski's theorem, one can see that for any ideal class $[I] \in \cl (\mathcal{O} _ d)$ one has $\operatorname{\frak{N}} ([I]) \ll \absolute {d} ^ {1/2}$ (except for the lower bound on the size of $\cl (\mathcal{O} _ d)$, all of these bounds are elementary and effective).

In \cite{Ellenberg-Michel-Venkatesh}, Ellenberg, Michel and Venkatesh prove Conjecture~\ref{joining conjecture} for $d _ i \to {-\infty}$ and $q _ i \in V_ {d} ^{*} (\Z) / \GL (2, \Z)$ with $\operatorname{\frak{N}} ([q_i]) \to \infty$ assuming one split prime  \emph{as long as $\operatorname{\frak{N}} ([q_i]) < \absolute {d _ i} ^ {1/2-\epsilon}$ for some fixed $\epsilon > 0$}. Essentially, their proof employed a variant of Linnik's method with a rather quantitative variant of Step~\ref{third step for Linnik} on p.~\pageref{third step for Linnik}, that gave first an equidistribution statement for the projection of an appropriate shift of the adelic torus $g_i [\T _ i (\mathbb A)]$ (notations as in the proof of Theorem~\ref{class group joining theorem}) to the quotient $\PGL (2, \R) / \Gamma _ i$ with $\Gamma _ i < \PGL (2, \Z)$ an appropriate congruence subgroup with $[\Gamma _ i: \PGL (2, \Z)] \to \infty$. This equidistribution of a single orbit in  a homogeneous space of increasing volume can then can be coupled with the equidistribution of the natural embedding 
\[\PGL (2, \R) / \Gamma _ i \hookrightarrow \left (\PGL (2, \R) / \PGL (2, \Z) \right) ^2\]
 as the uniform measure on a closed orbit\footnote{In fact, the graph of a Hecke correspondence.} of a diagonally embedded $\PGL(2,\R) \hookrightarrow \PGL(2,\R)^2$ to the uniform measure on the product space.

Assuming two split primes, Khayutin has been able to (essentially) prove Conjecture~\ref{joining conjecture} using a combination of ergodic and analytic tools, in particular Theorem~\ref{higher rank joining}:

\begin{Theorem}[Khayutin {\cite[Thm.~1.3]{Khayutin-CM-joining}}]\label{Khayutin theorem} Let $d _ i \to {-\infty}$ be a sequence of fundamental discriminants so that there are two primes $p,q$ for which $\bigl (\frac {d_i }{p}\bigr) = \bigl (\frac {d_i }{q}\bigr)= 1$ for all $i$, and let $q_i \in V ^{*} _ {d _ i} (\Z)$ satisfy $\operatorname{\frak{N}} ([q _ i]) \to \infty$.
Assume furthermore that the Dedekind $\zeta$-function of the fields $\Q(\sqrt{d_i})$ have no exceptional Landau-Siegel zero. Then the projection of the collection $\tilde V_{d_i}^{{\rm join},q_i}$
to $Y \times Y$ equidistribute as $i \to \infty$.
\end{Theorem}

It is very widely believed that Landau-Siegel zeros do no exist, as this is a very special, albeit important, case of the Riemann Hypothesis for Dedekind $\zeta$-functions. Moreover, even if these notorious Landau-Siegel zeros were to exist, they would have to be exceedingly rare: by a theorem of Landau (cf.~\cite [Thm.~5.28]{Iwaniec-Kowalski-book}, for any $A>1$ and $D$ large there will be at most one fundamental discriminant $d$ between $- D$ and $- D ^ A$ for which the Dedekind $\zeta$-function of $\Q(\sqrt{d})$ has such a zero.

Theorem~\ref{higher rank joining} allows one to deduce from Theorem~\ref{Khayutin theorem} the following seemingly much more general theorem:

\begin{Theorem}[Khayutin {\cite[Thm.~3.9]{Khayutin-CM-joining}}] Fix $k\in\N$. Let $d _ i \to {-\infty}$ be a sequence of fundamental discriminants so that there are two primes $p,q$ for which $\bigl (\frac {d_i }{p}\bigr) = \bigl (\frac {d_i }{q}\bigr)= 1$ for all $i$, and let $q_{i,1}, \dots, q_{i,k}\in V ^{*} _ {d _ i} (\Z)$ so that for any $1\leq j<\ell \leq k$ we have that $\operatorname{\frak{N}} ([q _ {i,j}]\odot [q _ {i,\ell}]^{-1}) \to \infty$.
Assume furthermore that the Dedekind $\zeta$-function of the fields $\Q(\sqrt{d_i})$ have no exceptional Landau-Siegel zero. Then the projection of the collections 
\[\left\{ ([q_{i,1}]\odot [q], \dots, [q_{i,k}]\odot [q]: q \in V _ {d_i} ^{*} (\Z) \right\}\]
to $Y^k$ equidistribute as $i \to \infty$.
\end{Theorem}

We end this rather long subsection with a striking equidistribution result by Aka, Einsiedler, and Wieser on the space $\grass _ {2, 4} (\R) \times Y ^ 4$ where $\grass _ {2, 4} (\R)$ is the projective variety (known as the \emph{Grassmanian}) of two-dimensional subspaces of a four dimensional space over $\R$, arising due to the ``accidental'' local isomorphism between $\SO (4)$ and $\SU (2) \times \SU (2)$ (or equivalently $\SO(3)\times\SO(3)$).

Consider now the quaternary\footnote{I.e. of 4 variables.} integer quadratic form $Q(x,y,z,w)= x ^2 + y ^2 + z ^2 + w ^2$, and consider all the binary integral quadratic forms $q$ of discriminant $d$ that can be represented by $Q$, i.e. so that there is a $4 \times 2$ integer matrix $M$ so that $q(x,y) \equiv Q \left (M \begin{pmatrix} x\\y \end{pmatrix} \right)$. Necessarily these binary quadratic forms will be positive definite, hence $d < 0$. It is a classical theorem that there are such binary quadratic forms iff $- d {} \not \equiv 0, 7, 12, 15 \bmod 16$. The image of $M$ is a rational subspace $L<\R^4$ of dimension two, hence in particular gives us a point in $\grass _ {2, 4}$.
The space perpendicular to $L$ also intersects the lattice $\Z ^ 4$ in a lattice, hence after choosing (arbitrarily) a basis for $L ^ \perp \cap \Z ^ 4$ we obtain another binary quadratic forms that can be shown to have also discriminant $d$.
Thus we obtain for any $d<0$, $- d {} \not \equiv 0, 7, 12, 15 \bmod 16$ a collection of triplets of points in $\grass _ {2, 4} \times (V _ d (\Z) / \GL (2, \Z)) ^2$, that we can project to $\grass _ {2, 4} \times Y ^2$. For each choice of a binary form $q$ represented by $Q$, Aka, Einsiedler, and Wieser magically pull two more rabbits (actually, only points of $Y$) out of the hat using the \emph{Klein map}, that assigns to any $L \in \grass _ {2, 4}$ a point in $(\mathbb S \times \mathbb S) / \{\pm 1\}$, with the $\pm 1$ acts by scalar multiplication on both factors. Identifying $\R ^ 4$ with the Hamilton quaternions $\left\{ x+iy +jz + ij w: x,y,z,w \in \R \right\}$, where $i^2=j^2=-1$, $ij=-ji$.
Pick two linearly independent vectors $v, w \in L$, and define
\begin{equation*}
a = v \overline w - \trace (v  \overline w) \qquad a' = w \overline v - \trace (w  \overline v)
.\end{equation*}
$a,a'$ are two traceless quaternions, hence lie in a three-dimensional space, and have the same norm. After rescaling, they define a point in $(\mathbb S \times \mathbb S) / \{\pm 1\}$ that turns out to be independent of the choice of generators $v$, $w$. If we chose $v,w$ to be generators of $L \cap \Z ^ 4$, then $a,a'$ will be integral and $Q(a)=Q(a')$ will be equal to $\absolute d$. The vectors $a$ and $a'$ give a continuous parameterization of $L$, so considering the joint distribution in the limit of the rescaled collections of triplets consisting of
\begin{itemize}
\item
a binary integral quadratic form $q$ of discriminant $d$ represented by $Q$,
\item the corresponding $L \in \grass _ {2, 4}$, and 
\item the quadratic form induced on $L^\perp$ 
\end{itemize}
already implicitly describes also the distribution of $a, a '$. However, these are integral vectors in a three-dimensional space, and the quadratic form induced by choosing generators of the lattices perpendicular to $a$ and $a'$ in this space gives as the desired two additional points of $Y$.

\begin{Theorem} [Aka, Einsiedler,Wieser \cite{Aka-Einsiedler-Wieser}] \label{Grassmanian theorem}
Let $d_i \to -\infty$ along the sequence $- d {} \not \equiv 0, 7, 12, 15 \bmod 16$ and so that $\bigl (\frac {d_i }{p}\bigr) = \bigl (\frac {d_i }{q}\bigr)= 1$ for two odd primes $p,q$. Then the collections of fivetuples in $\grass _ {2,4} \times Y ^ 4$ defined above become equidistributed as $i \to \infty$ with respect to the $\SO (4, \R) \times \PGL (2, \R) ^ 4$-invariant measure on this space.
\end{Theorem}

In a similar way to how Theorem~\ref{class group joining theorem} relates to Theorem~\ref{Aka-Einsiedler-Shapira theorem}, Theorem~\ref{Grassmanian theorem} is related to the distribution of the 6-tuples
\begin{multline*}
\Bigl\{ ([q],[q'], [q_{i,1}]\odot [q]\odot[q'], [q_{i,2}]\odot [q]\odot [q']^{-1},[q_{i,3}]\odot [q]^2,[q_{i,4}]\odot [q']^2): \\
q \in V _ {d_i} ^{*} (\Z) \Bigr\},
\end{multline*}
with $q_{i,1},\dots,q_{i,4} \in V _ {d_i} ^{*} (\Z)$ arbitrary, in $Y^6$.  We refer the reader to \cite[\S7]{Aka-Einsiedler-Wieser} for more details. Similarly to Theorem~\ref{class group joining theorem} and Theorem~\ref{Aka-Einsiedler-Shapira theorem}, the joining theorem, Theorem~\ref{higher rank joining}, is a key ingredient.

\subsection{Linnik's problem in $\PGL(3)$ and beyond}
Equidistribution of integer points discussed in \S\ref{Linnik's method} and \S\ref{beyond Linnik section} are at the core questions about the distribution of adelic points of $\Q$-tori on arithmetic quotients of forms of $\PGL (2)$, and their joinings. In this section we consider the more general question of density and equidistribution of ``adelic torus subsets'' --- sets of the form $g[\T (\mathbb A_F)]$ in $\mathbb G(\mathbb A_F)/\mathbb G (F)$ where $F$ is a number field, $\mathbb A_F$ is the Adele ring of $F$ and $\mathbb G$ is a reductive group over $F$, and $\T$ is an anisotropic $F$-torus. When $F=\Q$ we write $\mathbb A$ for $\mathbb A_\Q$.

Let $S$ be a set of places for $F$ (equivalence classes of embeddings of $F$ in local fields -- if $v$ is such an embedding, we denote by $F _ v$ the corresponding local field; we implicitly assume $F$ is dense in $F _ v$) containing at least one place $v$ in which $\T$ splits and all infinite places (embeddings of $F$ in $\R$ or $\C$ up to identifying conjugate embeddings in $\C$). Let $G = \prod_ {v \in S} \mathbb G (F _ v)$.

In \S\ref{beyond Linnik section}  essential (though mostly implicit) use was made of the fact that there is a natural projection from $\PGL (2, \mathbb A) / \PGL (2, \Q)$ to $\PGL (2, \R) / \PGL (\Z)$ with compact fibers. At the level of generality we are discussing here the picture is slightly more complicated: if one  would like to project $\mathbb G (\mathbb A _ F) / \mathbb G (F)$ to some quotient of $G$, in general one needs to take a finite number of lattices $\Gamma _ 1, \dots, \Gamma _ k  < G$ with $\Gamma _ i$ all conjugate over $\mathbb G (\Q)$ and all commensurable to $\mathbb G (\mathcal{O} _ {F, S})$ to obtain a natural projection
\begin{equation}\label{pi_S equation}
\pi _ S: \mathbb G (\mathbb A _ F) / \mathbb G (F) \to \bigsqcup ^{i = 1} _ k G / \Gamma _ i
\end{equation}
analogous to the projection from $\PGL (2, \mathbb A) / \PGL (2, \Q)$ to $\PGL (2, \R) / \PGL (\Z)$ (cf.~\cite[\S5]{Platonov-Rapinchuk}).\footnote{We already saw this phenomenon implicitly when discussing orthogonal groups in p.~\pageref{genus page} --- this is precisely the reason why for general ternary definite quadratic forms we need to consider not one quadratic form individually but the whole genus of quadratic forms locally equivalent to it.} We will also use $\pi _ S$ to denote the natural projection $\mathbb G (\mathbb A _ F)  \to  G $.

Since $\T$ was assumed to be $F$ anisotropic, the orbit $[\T (\mathbb A _ F)]$ in $\mathbb G (\mathbb A _ F) / \mathbb G (F)$ supports a $\T (A _ F)$-invariant probability measure $m _ {\T (\mathbb A _ F)}$.
The projection of $g [\T(\mathbb A _ f)]$ to $\bigsqcup ^{i = 1} _ k G / \Gamma _ i$ is a finite union of periodic $A_{g,S}=\pi_S(g) \prod_ {v \in S} \T (\Q _ S)\pi_S(g^{-1})$-orbits, and if $v \in S$ is a place where $\T$ is split the uniform measure on $\pi _ S (g[\T (\mathbb A _ F)])$, which is simply the average of the periodic measure on each of the $A_{g,S}$-periodic orbits comprising $\pi_S(g [\T(\mathbb A _ f)])$,  is invariant under a nontrivial $\Q _ v$-diagonalizable group.

Moreover if either the $\Q _ v$-rank of $\T (\Q _ v)$ is $\geq 2$, or $\T$ splits over at least one other place $v ' \in S$, $\pi_S(g [\T(\mathbb A _ f)])$ is invariant under a higher rank action, and one can hope to use tools of \S\ref{measure rigidity section} to study this set as well as the corresponding probability measure.

As an explicit example, we show how periodic $A$-orbits in $\PGL (n, \R) / \PGL (n, \Z)$ fit in this framework, where $A< \PGL(n, \R)$ is the full diagonal group. This corresponds to the above for $S=\{\infty\}$ and $\mathbb G = \PGL (n)$ when the shifting element is chosen appropriately. Indeed, an $A$-periodic orbit $A.[g_\infty]$ in $\PGL (n, \R) / \PGL (n, \Z)$ defines a $\R$-split $\Q$-tori in $\mathbb G$ by
\begin{equation*}
\mathbb T = C_{\mathbb G} \left (\PGL (n, \Z) \cap g_\infty ^{-1} A g_\infty \right)
.\end{equation*}
Then
\begin{equation}\label{adelic toral subset projection equation}
A.[g_\infty] \subseteq  \pi(g [\mathbb T (\mathbb{A})])\qquad\left(\text{$ g = (g_\infty, e, e,  \dots) \in \PGL (n, \mathbb{A})$}\right)
;\end{equation}
equality does not always occur in \eqref{adelic toral subset projection equation} --- in general, the projection $ \pi (g [\mathbb T  (\mathbb{A})])$ of the above adelic toral subset consists of a \emph{packet} of several periodic $A$-orbits $A.[g_\infty^{(i)}]$,  all with the same ``shape'', i.e. with the same
\[
\stab _ A ([g_\infty^{(i)}])= \left\{ a \in A: a . [g_\infty^{(i)}] = [g_\infty^{(i)}] \right\}
.\]

In this context, the following seems to be a natural conjecture. Conjecture~\ref{joining conjecture} can be viewed as a special case of this conjecture.

\begin{Conjecture}\label{adelic conjecture} Let $\mathbb G$ be a semi-simple algebraic group over a number field $F$, let $\mathbb \T _ i < \mathbb G$ be anisotropic $F$-tori, and let $g _ i \in \mathbb G (\mathbb A _ F)$. Let $\tilde {\mathbb G}$ be the simply connected cover of $\mathbb G$ with $j: \tilde {\mathbb G }\to \mathbb G$ the corresponding isogeny.
Then either:
\begin{enumerate}
\item \label{essentially equidistributed}
Any weak$^{*}$ limit of the uniform measures on $g _ i [\T (\mathbb A _ F)]$ is invariant under $j(\tilde {\mathbb G }(\mathbb A _ F))$.
\item \label{essentially restricted}
There exist a bounded sequence $h _ i \in \mathbb G (\mathbb A _ F)$ and a proper $\Q$-subgroup $\mathbb H < \mathbb G$ so that for infinitely many $i$,
\begin{equation*}
g _ i [\T (\mathbb A _ F)] \subset h _ i [\mathbb H (\mathbb A_ F)]
.\end{equation*}
\end{enumerate}
\end{Conjecture}

We remark that if $\mathbb G$ is simply connected, alternative~\eqref{essentially equidistributed} above is equivalent to $g _ i [\T (\mathbb A _ F)]$ become equidistributed in $\mathbb G (\mathbb A _ F) / \mathbb G (F)$, and in general implies that any weak$^*$-limit of the corresponding measures is homogeneous.
We also note that the assumption that $h _ i$ be bounded in alternative~\eqref{essentially restricted} of Conjecture~\ref{adelic conjecture} is equivalent to the following: there is a finite set of places $S$ so that for any $v \in S$ the $\mathbb G (F _ v)$-component of $h _ i$ remains bounded, and for any $v \not \in S$ the $\mathbb G (F _ v)$-component of $h _ i$ is in $\mathbb G (\mathcal O _ v)$, where $\mathcal O _ v$ is the maximal compact subring of $F _ v$.

Of interest are also results not just for the full group of adelic points $\T (\mathbb A _ F)$ but also large subgroups. These occur naturally in particular in the context of the study of special points on Shimura varieties: the orbit under the absolute Galois group on a special point, considered as a $\overline \Q$-point in an arithmetic model of $K\backslash G/\Gamma$, turns out to be such a group, though it is rather difficult to put ones hand on how big this group is --- cf. e.g.\ Tsimerman's paper \cite{Tsimerman-A-g} which proves an important special case of the Andr\'e-Oort conjecture (namely, when $\mathbb G=\operatorname{SP}(n)$) via such an analysis.
\medskip

The analogue of Conjecture~\ref{adelic conjecture} when $\mathbb T _ i$ is not a torus, e.g. when it is a semisimple or reductive group, is also highly interesting.\footnote{In general, in the context of $\Q$-groups or more generally groups defined over a number field~$F$, we will reserve $\T$ to denote an algebraic torus; we make an exception to this convention in this paragraph in order to abuse the notations of Conjecture~\ref{adelic conjecture} to cover a wider context.} If one assumes (implicitly or explicitly) that there is a fixed place $v$ for which $\mathbb T _ i (F_v)$ contains a unipotent subgroup  one can bring to bear a deep tools on unipotent flows, in particular Ratner's measure classification theorem \cite{Ratner-Annals} and its $S$-arithmetic generalizations by Ratner \cite{Ratner-padic} and by Margulis-Tomanov \cite{Margulis-Tomanov}, as was done by Eskin and Oh in \cite{Eskin-Oh-varieties}. Indeed, under some additional assumptions, and for $\T _ i$ semisimple, Einsiedler, Margulis and Venkatesh \cite{Einsiedler-Margulis-Venkatesh}, and these three authors jointly with Mohammadi \cite{Einsiedler-Margulis-Mohammadi-Venkatesh}, were able to give a \emph{quantitative} equidistribution statement; an exciting feature of \cite{Einsiedler-Margulis-Mohammadi-Venkatesh} is that thanks to the quantitative nature of the proof, it is even able to handle sequences of $\Q$-groups $\mathbb T _ i$ for which there is no place $v$ at which all (or even infinitely many) of these groups split. The discussion of these interesting works is unfortunately beyond the scope of this survey.

\medskip
Results towards Conjecture~\ref{adelic conjecture} for $\mathbb G = \PGL (n)$ or inner forms of $\PGL (n)$ (see bellow)  were obtained by Einsiedler, Michel, Venkatesh and the author in \cite{Einsiedler-Lindenstrauss-Michel-Venkatesh, Einsiedler-Lindenstrauss-Michel-Venkatesh-III} and strengthened in certain respects by Khayutin in \cite{Khayutin-double-torus}.

When considering periodic orbits of a fixed group $H$ on a space $X$, one can fix a Haar measure on $H$ and once this is done consistently measure the volume of all $H$ periodic orbits. When one allows the acting group to vary, one needs a slightly more sophisticated notion of volume:

\begin{Definition} [{\cite[Def.~4.3]{Einsiedler-Lindenstrauss-Michel-Venkatesh-III}}]
Let $\mathbb G$ be a fixed group defined over a number field $F$ and $\Omega \subset \mathbb G (\mathbb A _ F)$ a fixed neighborhood of the identity. Let $\mathbb H < \mathbb G$ be an $F$-subgroup.  We define the \emph{size}\footnote{In \cite{Einsiedler-Lindenstrauss-Michel-Venkatesh-III} we used the term ``volume'' of a periodic orbit to denote what we call here ``size''. We have decided to use a different terminology in this survey so we can unambiguously use ``volume'' to denote the volume of a periodic orbit with respect to a fixed Haar measure on the acting group.}  of an adelic shifted orbit  $g_i [\mathbb H (\mathbb A _ F)]$ to be $\infty$ if $\mathbb H (\mathbb A _ F) / \mathbb H (F)$ does not have finite $\mathbb H (\mathbb A _ F)$-invariant measure, and
\[
\size (g _ i [\mathbb H (\mathbb A _ F)]) = \frac { m_{\mathbb H (\mathbb A _ F)} (\mathbb H (\mathbb A _ F) / \mathbb H (F))}{m_{\mathbb H (\mathbb A _ F)} (\mathbb H (\mathbb A _ F) \cap g _ i ^{-1} \Omega g _ i) }
.\]
\end{Definition}

Note that changing $\Omega$ changes the size only up to a constant factor.

\begin{Theorem} [Einsiedler, Michel, Venkatesh and L. \cite{Einsiedler-Lindenstrauss-Michel-Venkatesh-III}]\label{ELMV3 theorem} Let ${\mathbb G}=\PGL (3)$, \ $F$ a number field, $g_i \in {\mathbb G} (\mathbb A _ F)$, and $\T _ i < {\mathbb G}$ a maximal $F$-torus. Assume that
\begin{enumerate}
\item\label{split condition ELMV} there is a place $v$ of $F$ so that (i)~$\T _ i$ is split over $F _ v$ and (ii)~$F _ v$ has no proper closed subfield --- i.e. either $F _ v \cong \R$ or $F _ v \cong \Q _ p$ for some prime $p$.\footnote {Part (ii) of this assumption was omitted in \cite{Einsiedler-Lindenstrauss-Michel-Venkatesh-III}, but is implicitly used in the proof. It is of course automatically satisfied for $F= \Q$.}
\item $\size (g_i [\T(\mathbb A _ F)]) \to \infty$ \label{volume going to infinity in ELMV}
\end{enumerate}
Then any limiting measure of the uniform distribution on $g_i [\T(\mathbb A _ F)]$ is invariant under the image of $\SL (3, \mathbb{A} _ F)$ in ${\mathbb G} (\mathbb{A} _ F)$.
In particular, if the class number of the integer ring in $F$ satisfies $\#\cl(\mathcal O_F)=1$ then the adelic torus sets $g_i [\T(\mathbb A _ F)]$ become equidistributed in ${\mathbb G} (\mathbb{A} _ F) / {\mathbb G} (F)$.
\end{Theorem}

Assumption \eqref{volume going to infinity in ELMV} in Theorem~\ref{ELMV3 theorem} turns out to be  equivalent to assumption \eqref{essentially restricted} in Conjecture~\ref{adelic conjecture}, as there are no proper $F$-subgroups of $\PGL(3)$ containing a maximal $F$-torus other than the torus itself, and the assumption that the volumes $\size (g_i [\T(\mathbb A _ F)]) \to \infty$ rules out the adelic torus subsets $g_i [\T(\mathbb A _ F)]$ being all in bounded translate of the image of a fixed $F$ torus in ${\mathbb G} (\mathbb{A} _ F) / {\mathbb G} (F)$.

The following is an easier to digest (weaker) form of Theorem~\ref{ELMV3 theorem} for $F = \Q$ where the adeles are not explicitly mentioned:

\begin{Theorem} [Einsiedler, Michel, Venkatesh and L. \cite{Einsiedler-Lindenstrauss-Michel-Venkatesh-III}]\label{ELMV3 simplified}
Let $A$ be the maximal diagonal group in $\PGL (3, \R)$. Let $V _ i$ be the sequence of all possible volumes of $A$-periodic orbits in $\PGL (3, \R) / \PGL (3, \Z)$ with respect to the Haar measure on $A$. For every $i$, let $\mathcal C _ i$ be the collection of $A$-periodic orbits of $A$ of volume exactly $V _ i$. Then these collections become equidistributed in $\PGL (3, \R) / \PGL (3, \Z)$ as $i \to \infty$.
\end{Theorem}

\noindent
The proof of Theorem~\ref{ELMV3 theorem} goes via a combination of analytic and ergodic tools: 
\begin{enumerate}[label={\arabic*.},ref={\arabic*}]
\item
Using the analytic theory of automorphic forms, specifically the subconvexity estimates of Duke, Friedlander, and Iwaniec \cite{Duke-Friedlander-Iwaniec-subconvexity-for-Artin} (or when $F$ is a general number field, an extension by Michel and Venkatesh of these subconvex bounds \cite{Michel-Venkatesh}), which incidentally are closely related to the works of Duke and Iwaniec mentioned in \S\ref{Linnik's method}, one shows that for some rather special functions $f \in L ^2 (\mathbb G (\mathbb{A} _ F) / \mathbb G (F))$
\begin{equation} \label{equidistribution for special functions}
\fint _ {g_i[\T _ i (\mathbb{A} _ F)]} f \to 0.
\end{equation}
Indeed, here the estimates are even quantitative.
\item Consider the $F_ v$-split tori $g_{i,v} \T _ i (\Q _ v) g _ {i, v} ^{-1}$, with $g_{i,v}$ denoting the $F_v$ component of $G_i$, and $v$ a place as in Theorem~\ref{ELMV3 theorem}.\eqref{split condition ELMV}. Without loss of generality these would converge to some $F_ v$-group $A_v$. If this group contains unipotent elements we can use Ratner's measure classification theorem (or more precisely its $S$-arithmetic extensions \cite{Ratner-padic, Margulis-Tomanov}). Otherwise one can use \eqref{equidistribution for special functions}, established for a certain collection of special $f$, to ensure that \emph{every ergodic component of any weak$^{*}$-limit} of the probability measures attached to $g_i[\T _ i (\mathbb{A} _ F)]$ has to have positive entropy with respect to the action of $A_v$, whence one can use the measure classification results of \cite{Einsiedler-Katok-Lindenstrauss, Einsiedler-Lindenstrauss-split}, e.g.~Theorem~\ref{split theorem simple}, to conclude the theorem.
\end{enumerate}

We note that Theorem~\ref{ELMV3 theorem} can be combined with the joining classification theorem Theorem~\ref{higher rank joining} to obtain joint equidistribution statements --- see \cite[Thm.~1.8]{Einsiedler-Lindenstrauss-joinings-2} for a precise statement.

\medskip

In \cite{Einsiedler-Lindenstrauss-Michel-Venkatesh}, a purely ergodic theoretic approach was used. This approach is not powerful enough to give a full equidistribution result, but on the other hand is significantly more flexible, and in particular gives information also about rather small subsets of an adelic torus subset. For simplicity (and to be more compatible with the terminology in \cite{Einsiedler-Lindenstrauss-Michel-Venkatesh}), we work over $\Q$ (instead of a general number field $F$) and use the more classical language of $A$-periodic orbits employed in Theorem~\ref{ELMV3 simplified}. We also assume for notational simplicity that the place where the $\Q$-tori we will consider is split is $\infty$, though the discussion below with minimum modification also holds for tori split over $\Q_p$ instead of $\R$.

Recall the relationship given in \eqref{adelic toral subset projection equation} between periodic $A$-orbits and the projection under $\pi _ S$ of appropriate shift of the adelic points of $\Q $-tori (with $\pi_S$ as in \eqref{pi_S equation} and $S=\{\infty\}$): a periodic orbit $A.[g]$ in $\PGL (n, \R) / \PGL (n, \Z)$ defines a $\Q$-torus $\T$ and $g.\pi _ S ([\T (\mathbb{A})])$ is a packet of periodic $A$-orbits of the same volume. In addition to the volume of a periodic $A$-orbit $A.[g]$, and its ``shape'' $ \stab _ A ([g])$, we can attach to this orbit an order in a totally real degree $n$-extension of $\R$, embedded in the subring of (not necessarily invertible) diagonal $n \times n$-matrices $D < M_{n\times n}(\R)$ as follows:
\begin{equation}
\label{order for periodic orbit equation}
\mathcal{O} _ {[g]} = \left\{ x \in D: x g \Z ^ n \subseteq g \Z ^ n \right\}
.
\end{equation}
The order $\mathcal{O} _ {[g]}$ in $D$ is best though of as an abstract order $\mathcal O$ in a totally real number field $K$ with $[K:\Q]=n$, together with an embedding $\tau$ of this order in $D$ (essentially this amounts to giving an ordering on the $n$ embeddings $K \to \R$).
The \emph{discriminant} $\disc (A.[g])$ of the periodic orbit $A.[g]$ is by definition the discriminant of the order $\mathcal{O} _ D$, i.e. (up to sign) the square of the co-volume of $\mathcal{O} _ {[g]}$ in $D$: it is an integer, since if $\alpha _ 1 = 1, \alpha _ 2, \dots, \alpha _ {n}$ are independent generators of $\mathcal{O} _ {[g]}$
\begin{equation*}
\disc (\mathcal{O} _ {[g]}) = \det (\trace  (\alpha _ i \alpha _ j))_{i,j=1}^n
.\end{equation*}
The relation between the volume of a periodic orbit $A.[g]$ (which is called by number theorists the \emph{regulator}) and the size of the adelic toral subset $\tilde {g}[\T (\mathbb{A})]$ (with $\tilde g$ the image of $g$ in $\PGL (n, \mathbb{A})$ under the obvious embedding $\PGL (n, \R) \hookrightarrow \PGL (n, \mathbb{A})$) is rather weak. Assuming the field generated by $\mathcal{O} _ {[g]}$ does not contain any nontrivial subfields\footnote{This assumption is needed only for the first inequality.}
\begin{equation}\label{sizes sizes sizes}
\log (\disc (A.[g])) ^ {n-1} \ll \vol (A.[g]) \ll \size \left (\bar {g}[\T (\mathbb{A})]\right ) = \disc (A.[g])^{1/2+o(1)}
\end{equation}
the last ``equality'' being ineffective (see \cite{Einsiedler-Lindenstrauss-Michel-Venkatesh, Einsiedler-Lindenstrauss-Michel-Venkatesh-III} for details).

Using the measure classification result in $\PGL (n, \R) / \PGL (n, \Z)$ of \cite{Einsiedler-Katok-Lindenstrauss} (a special case of Theorem~\ref{split theorem simple} above) and a rather crude entropy estimate the following was proved in \cite{Einsiedler-Lindenstrauss-Michel-Venkatesh}:

\begin{Theorem} [Einsiedler, Michel, Venkatesh and L. {\cite[Thm.~1.4]{Einsiedler-Lindenstrauss-Michel-Venkatesh}}] \label{ELMV GL theorem} Let $\Omega$ be a compact subset of $\PGL (n, \R) / \PGL (n, \Z)$, and $A$ the maximal diagonal subgroup of $\PGL (n, \R)$ for $n \geq 3$.
Then
\begin{equation*}
\#\left\{ \text{periodic $A$ orbits $A.[g] \subset \Omega$ with $\disc (A.[g]) \leq D$} \right\} \ll_{\epsilon, \Omega} D^ \epsilon .\end{equation*}
\end{Theorem}

Since the number of $A$-periodic orbits $A.[g]$ with $\disc (A.[g])$ is easily seen to be $\gg D^{c}$ for appropriate $c>0$, Theorem~\ref{ELMV GL theorem} can be viewed as evidence to the following conjecture, implied by Conjecture~\ref{C-SD conjecture}:

\begin{Conjecture}\label{periodic orbits conjecture} Let $n \geq 3$. Any compact $\Omega \subset \PGL (n, \R) / \PGL (n, \Z)$ contains only finitely many $A$-periodic orbits.
\end{Conjecture}

Conjecture~\ref{periodic orbits conjecture} follows from Conjecture~\ref{C-SD conjecture} using Cassels and Swinnerton-Dyer isolation result; cf.~\cite{Margulis-Oppenheim-conjecture}. Note that for $n=2$ the analogue of Theorem~\ref{ELMV GL theorem} is false; indeed, for any $\epsilon$ there is a compact $\Omega \subset \PGL (2, \R) / \PGL (2, \Z)$ so that
\begin{equation*}
\#\left\{ \text{periodic $A$ orbits $A.[g] \subset \Omega$ with $\disc (A.[g]) \leq D$} \right\} \gg D^ {1-\epsilon} ;\end{equation*}
cf.~\cite[Thm.~1.5]{Einsiedler-Lindenstrauss-Michel-Venkatesh}.

The non-compactness of $\PGL (n, \R) / \PGL (n, \Z)$ makes it harder to deduce a density statement from these rigidity results; however, for cocompact inner forms of $\PGL (n)$, namely $\PGL (1, \mathbb M)$ with $\mathbb M$ a central division algebra of degree $n$ over $\Q$, one can say more. Assume  $\mathbb M$ splits at $\R$, i.e. $\mathbb M \otimes \R \cong M _ {n \times n} (\R)$, and let $\mathcal O_{\mathbb M}$ be a maximal order in $\mathbb M(\Q)$. Then $\PGL (1, \mathbb M \otimes \R) \cong \PGL (n, \R)$ and $\PGL (1, \mathcal O_{\mathbb M})$ (or any subgroup of $\PGL (1, \mathbb M \otimes \R) $ commensurable to it) can be viewed as a cocompact lattice in $\PGL (n, \R)$.

Let $A$ be a maximal $\R$-split $\R$-torus in $\PGL (1, \mathbb M \otimes \R)$, and let $D$ be the abelian subalgebra of $\mathbb M \otimes \R$ commuting with $A$. As in \eqref{order for periodic orbit equation}, we can define for any periodic $A$-orbit $A .[g]$ in $\PGL (1, \mathbb M \otimes \R) / \PGL (1, \mathcal{O} _ {\mathbb M})$ an order in a totally real number field, and an embedding of this order to the algebra $D$ by considering
\begin{equation*}
\mathcal{O} _ {A.[g]} = D \cap g \mathcal{O} _ {\mathbf M} g ^{-1}
.\end{equation*}
While this will not be of relevance to our purposes, not all orders in totally real field of degree $n$ can appear in this way: indeed, an abstract orders $\mathcal{O} \cong \mathcal{O} _ {A.[g]}$ attached to periodic $A$-orbit in $\PGL (1, \mathbb M \otimes \R) / \PGL (1, \mathcal{O} _ {\mathbb M})$ has to satisfy the local compatibility condition that $\mathcal{O} \otimes \Q _ p$ can be embedded in $\mathbf M \otimes \Q _ p$ for all prime $p$.
In this context, one has the following:\footnote{The phrasing here is a bit stronger than that in \cite{Einsiedler-Lindenstrauss-Michel-Venkatesh}; the proof in \cite{Einsiedler-Lindenstrauss-Michel-Venkatesh} give this slightly stronger version.}

\begin{Theorem} [Einsiedler, Michel, Venkatesh and L. {\cite[Thm. 1.6]{Einsiedler-Lindenstrauss-Michel-Venkatesh}}] \label{ELMV compact theorem} 
Let~$\mathbb M$ be a division algebras over $\Q$ of degree $n$ so that $\mathbb M \otimes \R \cong M _ {n \times n} (\R)$.
Let $\mathcal{O} _ {\mathbb M}$ be a maximal order in $\mathbb M \otimes \Q$, and let $A$ be a maximal $\R$-split $\R$-torus in $\PGL (1, \mathbb M \otimes \R)$.
Let $\alpha > 0$, and for any $i$, let $\mathcal{C} _ i$ be a collection of (distinct) $A$-periodic orbits $\left\{ A.[g_{i,1}], \dots, A.[g_{i,k_i}] \right\}$ so that
\begin{equation*}
k_i \geq \left (\max_ {j} \disc (A \cdot [g _ {i, j}]) \right) ^ \alpha
\end{equation*}
and $k_i \to \infty$.
Assume that there is no subgroup $A \leq H < \PGL (1, \mathbb M \otimes \R)$ so that infinitely many $g_{i,j}$ lie on a single $H$-periodic orbit in $\PGL (1, \mathbb M \otimes \R) / \PGL (1, \mathcal{O} _ {\mathbb M})$. Then the collections
$\mathcal{C} _ i$ become dense in $ \PGL (1, \mathbb M \otimes \R) / \PGL (1, \mathcal{O} _ {\mathbb M})$, i.e.\ for every open $U \subset \PGL (1, \mathbb M \otimes \R) / \PGL (1, \mathcal{O} _ {\mathbb M})$ we have that there is an $i_0$ so that for $i>i_0$ there is a $j \in\{1,\dots,k_i\}$ so that $A.[g_{i,j}] \cap U\neq \emptyset$.
\end{Theorem}

To prove Theorem~\ref{ELMV compact theorem} one uses in addition to the ingredients used in Theorem~\ref{ELMV GL theorem}, namely Theorem~\ref{split theorem simple} and an appropriate entropy estimate, also a variant of the orbit closure/isolation theorems of  Weiss and the author \cite{Lindenstrauss-Barak}  and Tomanov \cite{Tomanov-maximal-tori} (cf.~Theorem~\ref{theorem with Barak}).
Applying Theorem~\ref{ELMV compact theorem} to the collections $\mathcal{C} _ i = \left\{ A.[g]: \mathcal{O} _ {A.[g]} = \tau_i(\mathcal{O} _ i) \right \}$ for $\mathcal{O} _ i$ a sequence of maximal orders in totally real degree $n$ number fields one gets the following, which can also be interpreted as a theorem about the projection of adelic toral subsets in $\PGL (1, \mathbb M \otimes \mathbb A) /\PGL (1, \mathbb M \otimes \Q)$ to $\PGL (1, \mathbb M \otimes \R) / \PGL (1, \mathcal{O} _ {\mathbb M})  $:

\begin{Corollary} In the notations of Theorem~\ref{ELMV compact theorem} let $\mathcal{O} _ i$ be the ring of integers in totally real fields $K _ i$, $\tau _ i$ embeddings of $\mathcal{O} _ i \hookrightarrow D$, and let
\begin{equation*}
\mathcal{C} _ i = \left\{ A.[g]: \mathcal{O} _ {A.[g]} = \tau _ i (\mathcal{O} _ i) \right\};
\end{equation*}
assume that $\mathcal{O} _ i$ are chosen so that the collections $\mathcal{C} _ i $ are nonempty. Assume that there is no fixed field $L$ of degree $d | n$ which is a subfield of infinitely many $K _ i$. Then the collections $\mathcal{C} _ i$ become dense in $\PGL (1, \mathbb M \otimes \R) / \PGL (1, \mathcal{O} _ {\mathbb M})  $.
\end{Corollary}

In \cite{Khayutin-double-torus}, a substantially more refined entropy estimate was given. This entropy estimates is quite interesting in its own sake, and in particular implies the following:

\begin{Theorem}[Khayutin \cite{Khayutin-double-torus}] Let $K_i$ be a sequence of totally real degree $n$ number fields and let $\mathcal{O} _ i$ be the ring of integers  $K _ i$. Let $\zeta$ be a generator for $K_i$ over $\Q$. Assume in addition that $n$ is prime and that the Galois group of the Galois extension of $K_i$  acts two-transitively on the Galois conjugates of $\zeta$. Let $\tau _ i$ be embeddings of $\mathcal{O} _ i \hookrightarrow D$,
\begin{equation*}
\mathcal{C} _ i = \left\{ A.[g]: \mathcal{O} _ {A.[g]} = \tau _ i (\mathcal{O} _ i) \right\};
\end{equation*}
and again assume that $\mathcal{O} _ i$ are chosen so that the collections $\mathcal{C} _ i $ are nonempty.
Then for any bounded continuous $f$ on $X=\PGL (1, \mathbb M \otimes \R) / \PGL (1, \mathcal{O} _ {\mathbb M})  $,
\begin{equation}\label{Khayutin estimate}
\liminf_ {i \to \infty} \frac {\sum_ {A.[g] \in \mathcal C_i} \int_ {A.[g]} f }{ \sum_ {A.[g] \in \mathcal C_i} \int_ {A.[g]} 1} \geq \frac 1{2(n-1)} \int_X f
.\end{equation}
\end{Theorem}

The techniques of \cite{Einsiedler-Lindenstrauss-Michel-Venkatesh} also imply an
estimate of the form \eqref{Khayutin estimate} but with a much worse bound. An important technical point is that the entropy bounds in \cite{Khayutin-double-torus}
also apply with regards to singular one parameter diagonal subgroups of $A$, hence would also be useful in the context of analyzing periodic orbits of a $\Q$-torus that is only partially split at a given place.

We remark that if one fixes a $\Q$-torus $\T$ and shifts it either in the real place or in one (or several) $p$-adic places one also obtain interesting equidistribution results, though they are now less related to diagonal flows and more to unipotent ones. See the work by Eskin, Mozes and Shah \cite{Eskin-Mozes-Shah-nondivergence,Eskin-Mozes-Shah-counting} for the former, with a nice application regarding counting matrices with a given characteristic polynomials in large balls in $\SL(n,\R)$, and \cite{Benoist-Oh} by Benoist and Oh for the latter, who use these results to study rational matrices with a given characteristic polynomial. Finally, we mention the work of Zamojski giving counting (and equidistribution) results for rational matrices in a given characteristic polynomial in terms of the height of these matrices \cite{Zamojski-thesis}. This leads to subtler issues, where unipotent flows or equidistribution of Hecke points do not apply. Instead, Zamojski uses measure rigidity of diagonal flows, building upon \cite{Einsiedler-Lindenstrauss-Michel-Venkatesh-III}. Notably, by fixing a $\Q$-torus, Zamojski is able to handle $\Q$-tori in $\SL(n,\R)$ for a general $n$; the fact that the $\Q$-torus is fixed allows one to avoid the need to use subcovexity results, and an additional averaging that is present in the problem studied by Zamojski allows handling intermediate subvarieties.

\section{Applications regarding quantum ergodicity}\label{sec: applications 2}

In this section, we consider applications of homogeneous dynamics, namely diagonal flows, to the study of Hecke-Maass cusp forms on $\H / \Gamma$ and their generalizations. We note that by the Selberg trace formula, Hecke-Maass forms can be considered as a dual object to the periodic $A$-trajectories considered in \S\ref{Linnik's method}, and though I am not aware of a dynamical result that makes use of this duality, the analogy is quite intriguing.

Consider first the case of $\mathbb G = \PGL (2)$. Then $K \backslash \mathbb G (\R) / \mathbb G (\Z)$ for $K = \operatorname {PSO} (2, \R)$ can be identified with the modular surface\footnote{To some, the modular curve\dots} $\H/\PSL (2, \Z)$.

To any primes $p$ there is a correspondence  --- the \emph{Hecke correspondence} which we will denote by $C ^ {\rm Hecke} _ p$   --- assigning to every $x \in  \H/\PSL (2, \Z)$ a set of $p+1$-points in this space. This correspondence can be described explicitly as follows: if $x = [z]$ for $z \in \H$ then
\begin{equation}\label{definition of Hecke correspondences}
C ^ {\rm Hecke} _ p ([z])= \Bigl\{ [pz],[z/p],[(z+1)/p], \dots,[(z+p-1)/p] \Bigr\}
;
\end{equation}
while each one of the points on the right hand side depends on the choice of representative $z$ of $[z]$ the collection of $p +1$-points is well-defined. Moreover this correspondence lifts to $\PGL (2, \R) / \PGL (2, \Z)$ giving to each $[x] \in \PGL (2, \R) / \PGL (2, \Z)$ a set (also denoted by $C ^ {\rm Hecke} _ p ([x])$) of $p+1$-points in $\PGL (2, \R) / \PGL (2, \Z)$ so that if $\pi _ Y:\PGL (2, \R) / \PGL (2, \Z) \to \H / \PSL (2, \Z)$ is the natural projection
\begin{equation*}
\pi _ Y \left (C ^ {\rm Hecke} _ p ([x]) \right) = C ^ {\rm Hecke} _ p (\pi _ Y ([x]))
.\end{equation*}
An important property of the Hecke correspondence is its reflexivity:
\begin{equation} \label{reflexivity equation} \text{$[y] \in C ^ {\rm Hecke} _ p ([x])$ iff $[x] \in C ^ {\rm Hecke} _ p ([y])$}
.\end{equation}
moreover $C ^ {\rm Hecke} _ p ([\bullet])$ is equivariant under left translations on $\PGL (2, \R) / \PGL (2, \Z)$, i.e. 
\[
C ^ {\rm Hecke} _ p (h.[x]) = h. C ^ {\rm Hecke} _ p ([x])
\]
which implies that on $\H / \PGL (2, \Z)$ each branch of $C ^ {\rm Hecke} _ p ([\bullet])$ is a local isometry.

In terms of the projection (for $S = \left\{ \infty \right\}$)
\[\pi _ S: \PGL (2, \mathbb{A}) / \PGL (2, \Q) \to \PGL (2, \R) / \PGL (2, \Z),\] 
 if $a_p \in \PGL (2, \mathbb{A})$ is the element equal to $\begin{pmatrix} p& \\& 1 \end{pmatrix}$ in the \mbox{$\Q _ p$-component} and the identity in every other component then for any~$[x] \in \PGL (2, \R) / \PGL (2, \Z)$
\begin{equation}\label{definition of Hecke operators using adels}
C ^ {\rm Hecke} _ p ([x]) = \pi _ S (a_p.\pi _ S ^{-1} ([x]))
.\end{equation}
Phrased slightly differently, if we consider an $a_p$ orbit
\[
 \{[\bar x], a_p.[\bar x], \dots, a_p^k.[\bar x]\} \subset  \PGL (2, \mathbb{A}) / \PGL (2, \Q)
\]
 and project it to $\PGL (2, \R) / \PGL (2, \Z)$ we will get a sequence of points $[x_0]$, \dots,$[x_k]$
with $[x_i] \in C ^ {\rm Hecke} _ p ([x_{i-1}])$; moreover it can be shown that this discrete trajectory is ``non-backtracking'' in the sense that $[x_i] \neq [x_{i+2}]$.

Using the Hecke correspondences $C ^ {\rm Hecke} _ p (\bullet)$ on $\H / \PSL (2, \Z)$ we define for any prime $p$ a self-adjoint operator $T _ p$, called \emph{Hecke operators}, on $L ^2 (\H / \PSL (2, \Z))$ by
\begin{equation*}
(T _ p f) (x) = p ^ {-1/2} \sum_ {y \in C ^ {\rm Hecke} _ p (x)} f (y)
.\end{equation*}
It follows from the relation between the Hecke correspondences and actions of diagonal elements in $\PGL (2, \mathbb{A}) / \PGL (2, \Q)$ (or directly from the definition of these correspondences in \eqref{definition of Hecke correspondences}) that for every prime $p$, $q$, the operators $T _ p$ and $T _ q$ commute. Moreover, using the symmetry of the Hecke correspondences \eqref{reflexivity equation} and the fact that each branch of $C ^ {\rm Hecke} _ p (\bullet)$ is a locally isometry, one sees that operators $T _ p$ are self adjoint operators on $L ^2 (\H / \PSL (2, \Z))$ commuting with the Laplacian. Thus using the well-known spectral properties of the Laplacian the discrete spectrum of $\Delta$ in $L ^2 (\H / \PSL (2, \Z))$ is spanned by joint eigenfunctions of $\Delta$ and all Hecke operators. Moreover except for finitely many of them (in fact only the constant function), these eigenfunctions will be \emph{cusp forms}, i.e. eigenfunctions of $\Delta$ with the property that their integral over any periodic horocycle on $\H / \PSL (2, \Z)$ is zero. These joint eigenfunctions are called \emph{Hecke-Maass cusp forms}.\footnote{For us cuspidality of the forms is not relevant - only that these are eigenfunctions of $\Delta$ and all $T_p$ in $L ^2 (\H / \PSL (2, \Z))$.}

A similar setup works also in cocompact quotients. Let $\mathbb M$ be a quaternionic division algebra over $\R$ with $\mathbb M \otimes \R \cong M _ {2 \times 2} (\R)$. Let $\mathbb G = \PGL (1, \mathbb M)$. Then $\mathbb G (\R) \cong \PGL (2, \R)$, and if $\mathcal{O} _ \mathbb M$ is a maximal order in $\mathbb M$ then $\Gamma = \mathcal{O} _ \mathbb M ^{\times} / \Q ^{\times}$ is a cocompact lattice in $\mathbb G (\R)$ commensurable to $\mathbb G (\Z)$.\footnote {To define  $\mathbb G (\Z)$ (at least in the way we do it in the survey) one needs to choose a $\Q$-embedding of $\mathbb G$ in some $\SL (N)$; reasonable people might do this in different ways, but they would all agree that the $\Gamma$ we defined is commensurable to $\mathbb G (\Z)$.}
Taking as before $S = \left\{ \infty \right\}$ then by \eqref{pi_S equation}
\begin{equation*}
\pi _ S: \mathbb G (\mathbb A) / \mathbb G (\Q) \to \bigsqcup _{i = 1} ^ k \PGL (2, \R)/ \Gamma _ i.
\end{equation*}
Indeed, maximal orders in $\R$-split quaternion algebras have class number 1, so in fact the image is a single quotient $\PGL (2, \R)/ \Gamma$, though if one takes a non-maximal order $\mathcal O$ a disjoint union is needed (cf.~\cite[\S2.2]{Rudnick-Sarnak} and references therein).
At any place $p$ in which $\mathbb M \otimes \Q _ p \cong M _ {2 \times 2} (\Q _ p)$ (in particular, for all but finitely many places) we can choose (non-canonically) an element $a _ p$ as in the paragraph above \eqref{definition of Hecke operators using adels}, and this allows us using \eqref{definition of Hecke operators using adels}  to define Hecke correspondences $C ^ {\rm Hecke} _ p (\bullet)$ on $\ \PGL (2, \R)/ \Gamma $ and $\H/ \Gamma$ as well as a family of self adjoint operators $T_p$ commuting with each other and with $\Delta$ on $L ^2 (\H/ \Gamma)$. Then $L ^2 ( \H/ \Gamma)$ is spanned by Hecke-Maass forms --- joint eigenfunctions of $\Delta$ and all $T _ p$.\footnote {In this case, these joint eigenfunctions, even the constant function, satisfy the condition of cuspidality automatically (if somewhat vacuously) since there are no periodic horospheres!}

Motivated in part by the study of Hecke-Maass forms, Rudnick and Sarnak made the following bold conjecture regarding any hyperbolic surface:
\begin{Conjecture} [Quantum Unique Ergodicity\footnotemark Conjecture \cite{Rudnick-Sarnak}] \label{quantum unique ergodicity conjecture}\footnotetext{Abbreviated QUE bellow.}
Let $M$ be a compact manifold of negative sectional curvature. Let $\left\{ \phi _ i \right\}$ be a complete orthonormal sequence of eigenfunctions of the Laplacian $\Delta$ on $M$ ordered by eigenvalue. Then the probability measures $\absolute {\phi _ i(x)}^2 \,d \vol _M(x)$ converge weak$^*$ to the uniform measure on $M$.
\end{Conjecture}

In their paper, Rudnick and Sarnak focus on the case of Hecke-Maass forms, showing that any weak$^*$ limit of a subsequence $\absolute {\phi _ {i_j}(x)}^2 \,d \vol _M(x)$ cannot be supported on finitely many closed geodesics. The multiplicities in the spectrum of the Laplacian on arithmetic surfaces $\H / \Gamma$ with $\Gamma$ as above (or more generally the $\bigsqcup _{i = 1} ^ k \H/ \Gamma _ i$ on which the Hecke correspondences are defined if we work with non-maximal orders) are not well understood.
 Empirically, these multiplicities seem to be bounded, indeed in favorable cases, the multiplicity of every eigenvalue of $\Delta$ seems to be one, so one does not seem to lose much by using a sequence of Hecke-Maass forms instead of an arbitrary sequence of eigenfunctions of $\Delta$. 
 Let $M$ be an arithmetic surface (or a more general local symmetric manifold $K \backslash G / \Gamma$ with $K<G$ maximal compact, $G= \mathbb G (\R)$ for $\mathbb G$ a semisimple $\Q$-group, and $\Gamma < G$ a congruence lattice [in particular commensurable to $\mathbb G (\Z)$]). 
 We shall call the closely related question to Conjecture~\ref{quantum unique ergodicity conjecture}, of whether on such an $M$ the probability distributions $\absolute {\phi _ i(x)}^2 \,d \vol _M(x)$ corresponding to Hecke-Maass forms $\phi _ i$  (i.e. joint eigenfunctions of $\Delta$ and all Hecke operators) converge weak$^*$ to the uniform measure the \emph{Arithmetic Quantum Unique Ergodicity Problem}.

Conjecture~\ref{quantum unique ergodicity conjecture} is to be compared to the following ``quantum ergodicity'' theorem of Schnirelman, Colin de Verdiere, and Zelditch:

\begin{Theorem} [Schnirelman, Colin de Verdiere, Zelditch \cite{Schnirelman,Colin-de-Verdiere-quantum,Zelditch-uniform-distribution}] Let $M$ be a compact manifold so that the geodesic flow on the unit tangent bundle of $M$ is ergodic. Let $\left\{ \phi _ i \right\}$ be a complete orthonormal sequence of eigenfunctions of the Laplacian $\Delta$ on $M$ ordered by eigenvalue. Then there is a subsequence $i _ j$ of density one so that restricted to this subsequence  $\absolute {\phi _ {i_j}(x)}^2 \,d \vol _M(x)$ converge weak$^*$ to the uniform measure on $M$.
\end{Theorem}

Using Theorem~\ref{one parameter plus recurrence rigidity theorem} as well as an entropy estimate by Bourgain and the author \cite{Bourgain-Lindenstrauss} the author has been able to prove the following, in particular establishing Arithmetic Quantum Unique Ergodicity for compact hyperbolic surfaces:

\begin{Theorem} [\cite{Lindenstrauss-quantum}]\label{arithmetic quantum unique ergodicity theorem} Let $\phi _ i$ be an $L^2$-normalized sequence of Hecke-Maass forms on an arithmetic surface\footnote{More generally, finite union of surfaces $\bigsqcup _{i = 1} ^ k \H/ \Gamma _ i$.} $M=\mathbb H / \Gamma$  with the lattice $\Gamma$ either a congruence subgroup of $\PGL (2, \Z)$ or arising from an order in a quaternion division algebra over $\Q$ as above. Suppose $\absolute {\phi _ {i}(x)}^2 \,d \vol _M(x)$ converge weak$^{*}$ to a measure~$\mu$ on $\mathbb H / \Gamma$. Then $\mu$ is up to a multiplicative constant the uniform measure on $\mathbb H / \Gamma$. In particular, arithmetic quantum unique ergodicity holds for compact arithmetic surfaces.
\end{Theorem}

What is not addressed in that theorem is the question whether there can be escape of mass for the sequence of measures $\absolute {\phi _ {i}(x)}^2 \,d \vol _M(x)$ for $\Gamma$ a congruence sublattice of $\PGL (2, \Z)$. What is shown by Theorem~\ref{arithmetic quantum unique ergodicity theorem} is that whatever remains converges to the uniform measure. This difficulty was resolved by Soundararajan using an elegant analytic argument:
\begin{Theorem} [Soundararajan \cite{Sound-QUE}]\label{Sound's theorem} Let $\phi$ be a Hecke-Maass form on $\H / \PSL (2, \Z)$, normalized to have $L ^2$-norm one. Then
\begin{equation*}
\int_ {\substack{\absolute{x} \leq 1/2\\y>T}} \absolute{\phi (x+iy)}^2 \,d \vol_{\H} (x+iy) \ll \frac {\log (e T) }{ \sqrt T}
\end{equation*}
with an absolute implicit constant.
\end{Theorem}

Theorem~\ref{Sound's theorem} implies in particular that any weak$^{*}$ limit of a sequence of measures corresponding to Hecke-Maass forms $\absolute {\phi _ {i}(x)}^2 \,d \vol _M(x)$ is a probability measure, hence using Theorem~\ref{arithmetic quantum unique ergodicity theorem} Arithmetic QUE holds also in the non-compact case where $\Gamma=\PSL (2, \Z)$ (a similar argument works for congruence sublattices of $PSL (2, \Z)$).

\medskip

The entropy bound by Bourgain and the author in \cite{Bourgain-Lindenstrauss} gives a uniform upper bound on measures of small balls in $\PGL(2,\R)/\Gamma$ of an appropriate lift of the measures $\absolute {\phi _ {i}(x)}^2 \,d \vol _M(x)$ to $\PGL(2,\R)/\Gamma$. An alternative approach by Brooks and the author using only one Hecke operator gives a less quantitative entropy statement that is still sufficient to prove quantum unique ergodicity:
\begin{Theorem} [Brooks and L. \cite{Brooks-Lindenstrauss-que}] Let $\phi _ i$ be an $L^2$-normalized sequence of smooth functions on $\H / \Gamma$ with $\Gamma$ an arithmetic co-compact lattice arising from an order in a quaternionic division algebra over $\Q$ as above. Assume that for some sequences $\lambda _ i \to \infty $, $\lambda _ {i,p} \in \R$, $\omega _ i \to 0$
\begin{align*}
\norm {\Delta \phi _ i - \lambda _ i \phi _ i}_2 \leq \lambda _ i ^ {1/2} \omega _ i \qquad \norm {T _ p \phi _ i - \lambda _ {i,p} \phi _ i}_2 \leq \omega _ i
.\end{align*}
Then $\absolute {\phi _ i(x)} ^2 \,d \vol (x)$ converge weak$^{*}$ to the uniform measure on $\H / \Gamma$.
\end{Theorem}

\medskip

A surprising link between quantum unique ergodicity and the number of nodal domains for Hecke-Maass forms $\phi$ on $\H / \Gamma$ was discovered by Jang and Jung. If $\phi: M \to \R$ is a $\Delta$-eigenfunction, say $\Delta \phi + \lambda \phi = 0$, on a compact surface $M$, Courant's Nodal Domain Theorem and the Weyl Law imply that the number of nodal domain $\mathcal N (\phi)$ for $\phi$ satisfies $N (\phi) \ll \lambda$. However it is well-known that in general $\mathcal N (\phi)$ could be much less: indeed in the two sphere $\mathbb S$ there is a sequence of $\Delta$-eigenfunctions with eigenvalues $\to\infty$ with $\mathcal N (\phi)\leq 3$, and in general it is very hard to bound the number of nodal domains from below; for more details, cf.\ e.g.~\cite{Jang-Jung} and the references given by that paper.

\begin{Theorem} [Jang and Jung \cite{Jang-Jung}] let $\phi _ i$ be a sequence of Hecke-Maass forms on $\H / \Gamma$ for $\Gamma$ an arithmetic triangle group. Then $\lim_ {i \to \infty} \mathcal N (\phi _ i) = \infty$.
\end{Theorem}

Triangle groups are discrete subgroups of the isometry group of $\H$ generated by reflections in three sides of a triangle with angles $\pi/a, \pi/b, \pi/c$.  To such a group $\Gamma '$ we can attach the orbifold $\H / \Gamma$ where $\Gamma <\Gamma '$ is the group of orientation preserving isometries. $\Gamma '$ is generated by $\Gamma$ and a reflection $\sigma$ with $\sigma \Gamma \sigma = \Gamma$, hence $\sigma$ induces a orientation reversing involution on $\H / \Gamma$ (for convenience we will also call $\Gamma$ a triangle group). An arithmetic triangle group (there are only finitely many of these) is a triangle group that is commensurable to $\PSL (2, \Z)$ or a lattice coming from a quaternionic order over $\Q$ as above. Examples are $\Gamma = \PSL (2, \Z)$ itself (giving the triangle group $(\infty,3,3)$) and the compact triangle group $(2,6,6)$.
Quantitative results giving $\mathcal N (\phi)\gg \lambda ^ {\frac 1 {27}-\epsilon}$ were given by Ghosh, Reznikov and Sarnak, \emph{assuming the Lindel\"of hypothesis for L-functions of $\GL (2)$-forms} \cite{Ghosh-Reznikov-Sarnak, Ghosh-Reznikov-Sarnak-II}; quantum unique ergodicity is used as a (partial) substitute to the Lindel\"of hypothesis in \cite{Jang-Jung}.

\medskip

In addition to considering a sequence of Hecke-Maass forms with eigenvalue $\lambda _ i \to \infty$, one can consider on a given quotient $\H / \Gamma$ a sequences of holomorphic Hecke forms of weight $\to \infty$. These correspond to irreducible $\PGL (2, \R)$-representation in $L ^2 (\PGL (2, \R) / \Gamma)$ which have no $\operatorname{PO} (2, \R)$-invariant vectors. Using  analytic techniques, and in particular making heavy use of the Fourier expansion of such forms in the cusp, Holowinsky and Soundararajan \cite{Holowinsky-Sound} were able to prove an arithematic quantum unique ergodicity theorem for these automorphic forms \emph{for $\Gamma$ a congruence subgroup of $\PGL (2, \Z)$.} The techniques of Holowinsky and Soundararajan seem to be restricted to the noncompact case; the analogous question for compact quotients, and even on the sphere $\mathbb S$, remain important open questions.

A related question involves fixing the weight (or bounding the Laplacian eigenvalue) but considering a sequence of ``newforms'' on a tower $\H / \Gamma (N)$ of congruence subgroups.
This question makes sense also in a purely discrete setting: instead of taking a quaternion division algebra $\mathbb M$ over $\Q$ which is split at infinity, one can take a definite quaternion algebra such as the Hamilton quaternions 
\[\mathbb M = \Q + i \Q + j \Q + i j \Q \qquad i ^2 = j ^2 = -1, ij=-ji.\]
 For such $\mathbb M$, the group $\PGL (1, \mathbb M \otimes \R)$ is compact (in fact, isomorphic to $\SO (3, \R)$), and so for $S= \left\{ \infty, p \right\}$, the $S$-arithmetic projection
\begin{equation}
\label{projection p infinity}
\pi _ S: \PGL (1, \mathbb M \otimes \mathbb A) / \PGL (1, \mathcal{O} _ {\mathbb M}) \to \PGL (1, \mathbb M \otimes \R) \times \PGL (1, \mathbb M \otimes \Q _ p) / \Gamma
\end{equation}
can be composed with a further projection by dividing the R.H.S.\ of~\eqref{projection p infinity} from the left by the compact group $\PGL (1, \mathbb M \otimes \R) \times K(n)$ with $K(n) \leq   \PGL (1, \mathbb M \otimes \Z _ p)$ a compact open subgroup.
This gives a map from $\PGL (1, \mathbb M \otimes \mathbb A) / \PGL (1, \mathcal{O} _ {\mathbb M})$ to a finite set. For every $q\neq p$ the $q$-Hecke correspondence gives this finite set the structure of a $q+1$-regular graph --- the Lubotzky, Phillips and Sarnak ``Ramanujan graphs'' \cite{Lubotzky-Phillips-Sarnak-Ramanujan-graphs}. Taking Hecke newforms on a sequence of these graphs with decreasing $K(n)$, Nelson \cite{Nelson-GL2p} was able to use an adaptation of the method of proof of Theorem~\ref{arithmetic quantum unique ergodicity theorem} to prove Arithmetic QUE in the level aspect for newforms corresponding to principle series representations of $\PGL(2,\Q_p)$. The restriction to principle series representations is the analogue in this context of the restriction in Theorem~\ref{arithmetic quantum unique ergodicity theorem} to Maass forms (i.e.~Laplacian eigunfunctions) as opposed to holomorphic forms.

The dynamical approach to arithmetic quantum unique ergodicity can be extended to other arithmetic quotients. Notable work in that direction was done by Silberman and Venkatesh:

\begin{Theorem} [Silberman and Venkatesh \cite{Silberman-Venkatesh-I, Silberman-Venkatesh-II}]\label{Silberman Venkatesh theorem} Let $\mathbb G = \PGL (1, \mathbb M)$ where $\mathbb M$ is a degree $n$-division algebra over $\Q$, split over $\R$, for $n \geq 3$. Let $G= \mathbb G (\R) \cong \PGL (n, \R)$, $K <G$ maximal compact, and $\Gamma$ a lattice in $G$ arising from a maximal order in $\mathbb M$. Let $\phi _ j$ be a sequence of $L^2$-eigenfunctions of the ring of invariant differential operators on $K \backslash G$ as well as all Hecke-operators. Assume that the irreducible $G$-representation $H_j < L ^2 (G/\Gamma)$ of $G$ spanned by left translations of $\phi _ j$ is a principle series representation with parameter $ \mathbf t_j \in i \mathfrak a _\R$ that stays away (uniformly in $j$) from the edges of the positive Weyl chamber in $i \mathfrak a _\R$. Then $\absolute {\phi _ j(x)}^2 d\vol_{K \backslash G/\Gamma}$ converge to the uniform measure on $\absolute {\phi _ j(x)}^2 d\vol_{K \backslash G/\Gamma}$.
\end{Theorem}

The proof of Theorem~\ref{Silberman Venkatesh theorem} proceeds, similarly to that of Theorem~\ref{arithmetic quantum unique ergodicity theorem}, by ``lifting'' the probability measures $\absolute {\phi _ j(x)}^2 d\vol_{K \backslash G/\Gamma}$ to an (approximately) $A$-invariant probability measure on $G / \Gamma$. This part of the argument, which is carried out in \cite{Silberman-Venkatesh-I}, uses only the information about the behavior of $\phi _ j$ at the infinite place. Then using the information about other places, explicitly the fact that $\phi _ i$ are eigenfunctions of all Hecke operators, positive entropy of any ergodic component of any weak$^*$ limit as above is derived in \cite{Silberman-Venkatesh-II}; this argument is related to the entropy estimate of Bourgain and the author in \cite{Bourgain-Lindenstrauss}. Once this entropy estimate is established, one can employ the measure classification results of \cite{Einsiedler-Katok-Lindenstrauss} (special case of Theorem~\ref{split theorem simple}) to deduce the above arithmetic quantum unique ergodicity result.

\section{Applications regarding Diophantine approximations}\label{sec: applications 3}

We started this survey with a historical introduction concerning some of the origins of the study the rigidity properties of higher rank diagonal actions. One important such work was the paper \cite{Cassels-Swinnerton-Dyer} of Cassels and Swinnerton-Dyer, relating Littlewood's Conjecture (Conjecture~\ref{Littlewood Conjecture}) to Conjecture~\ref{C-SD conjecture} regarding bounded $A$-orbits in $\PGL (3, \R) / \PGL (3, \Z)$.

It is therefore not surprising that the significant progress obtained towards understanding these higher rank diagonal actions in the half-century since Cassels and Swinnerton-Dyer's seminal paper shed some light on Diophantine questions, though Littlewood's Conjecture itself remains at present quite open. Indeed, in terms of concrete (e.g. algebraic) numbers $\alpha, \beta \in \R$ for which Littlewood's conjecture can be verified, i.e. so that 
\begin{equation}\label{Littlewood is satisfied}
\inf_ {n>0} n \norm {n \alpha} \norm {n \beta} = 0,
\end{equation}
I am not aware of any nontrivial examples beyond that given by Cassels and Swinnerton-Dyer in \cite{Cassels-Swinnerton-Dyer}, namely those $\alpha,\beta \in \R$ so that the field $\Q(\alpha,\beta)$ they generate is of degree 3  over $\Q$.

The following was proved by Einsiedler, Katok and the author in \cite{Einsiedler-Katok-Lindenstrauss} using measure rigidity of the action of the diagonal group on $\PGL (3, \R) / \PGL (3, \Z)$:

\begin{Theorem} [Einsiedler, Katok, L.]\label{theorem towards Littlewood} For any $\epsilon > 0$ the set
\begin{equation}\label{a subset of the square}
\left\{ (\alpha, \beta \in [0, 1] ^2:\inf_ {n>0} n \norm {n \alpha} \norm {n \beta} \geq \epsilon \right\}
\end{equation}
has zero (upper) box dimension.
\end{Theorem}

Zero upper box dimension simply means that the set~\eqref{a subset of the square} can be covered by $\ll_{\delta,\epsilon} N^{\delta}$ squares of diameter $N^{-1}$ for any $N>0$. This of course implies that the set of exceptions to Littlewood conjecture has Hausdorff dimension zero, and moreover for any $\alpha \in \R$ outside a set of Hausdorff dimension zero, \eqref{Littlewood is satisfied} holds for every $\beta \in \R$.

This latter statement can actually be made more explicit: let $\alpha \in [0,1]$ be given. Expand $\alpha$ to a continued fraction
\begin{equation}\label{continued fractions}
\alpha = \cfrac{1}{n_1+\cfrac{1}{n_2+\cfrac{1}{n_3+\cfrac{1}{n_4+\dots}}}}
.
\end{equation}
If the sequence $n _ d$ is unbounded then already $\inf_ {n>0} n \norm {n \alpha} = 0$ and hence \eqref{Littlewood is satisfied} holds for every $\beta$. For any $k$ let $N_k (\alpha)$ denotes the number of possible $k$-tuples of integers $i_1,\dots,i_k$ appearing in the continued fraction expansion of $\alpha$, i.e.\ so that there is some $\ell \in \N$ so that
\[
(i_1,\dots,i_k)=(n_{\ell},\dots, n _ {\ell + q -1})
.\]

The following proposition follows readily from the techniques of \cite{Einsiedler-Katok-Lindenstrauss}; we leave the details the imagination of the interested reader, but the key point is that the condition given in the proposition on $\alpha$ can be used to verify the positive entropy condition.

\begin{Proposition} Let $\alpha \in [0,1]$ be such that the continued fraction expansion of $\alpha$ satisfies that
\begin{equation*}
\lim_ {q \to \infty} \frac {\log (N _ k (\alpha))} {k} > 0
\end{equation*}
(this limits exists by subadditivity). Then for any $\beta \in \R$ equation \eqref{Littlewood is satisfied} holds.
\end{Proposition}

de Mathan and Teuli\'e gave the following analogue to Conjecture~\ref{Littlewood Conjecture}

\begin{Conjecture}[de Mathan and Teuli\'e \cite{de-Mathan-Teulie}] \label{p-adic Littlewood conjecture}
For any prime $p$ and any~$\alpha \in \R$
\begin{equation} \label{p-adic Littlewood equation}
\inf_ {n > 0} n \absolute {n} _ p\norm {n \alpha}  = 0.
\end{equation}
\end{Conjecture}

Recall that $\absolute {n} _ p = p^{-k}$ if $p^k \mid n$ but $p^{k+1} \nmid n$, hence \eqref{p-adic Littlewood equation} is equivalent to
\begin{equation*}
\inf_ {n > 0, k\geq 0} n\norm {n p^k \alpha} = 0
.\end{equation*}
Note that by Furstenberg's theorem (Theorem~\ref{Furstenberg's theorem})\footnote{Another result we cited in the introduction that played a central role in the development of the subject!} for any two distinct primes~$p,q$
\begin{equation*}
\inf_ {n > 0} n \absolute {n} _ p\absolute {n} _ q\norm {n \alpha}  = 0,
\end{equation*}
since either $\alpha$ is rational, in which case $\liminf \norm {n \alpha}=0$, or $\{p^kq^\ell \alpha \bmod 1\}$ is dense in $[0,1]$, in particular $\inf_{k,\ell\geq0} \norm {p^k q^\ell \alpha}=0$.

By a variant of the argument of Cassels and Swinnerton-Dyer, de Mathan and Teuli\'e show \eqref{p-adic Littlewood equation} holds for quadratic irrational $\alpha \in \R$. Interestingly, Adiceam, Nesharim and Lunnon give in \cite{Adiceam-Nesharim-Lunnon} a completely explicit (and non-obvious) counterexample to the function field analogue of Conjecture~\ref{p-adic Littlewood conjecture}, also stated in \cite{de-Mathan-Teulie}. Using a similar argument to \cite{Einsiedler-Katok-Lindenstrauss}, but using the measure classification result of \cite{Lindenstrauss-quantum} instead of that in \cite{Einsiedler-Katok-Lindenstrauss}, Einsiedler and Kleinbock prove in \cite{Einsiedler-Kleinbock} that for any $\epsilon > 0$ the set of $\alpha \in [0, 1]$ for which
\begin{equation*}
\inf_ {n > 0} n \absolute {n} _ p\norm {n \alpha} \geq \epsilon
\end{equation*}
has zero box dimension.

Theorem~\ref{theorem towards Littlewood}, unlike many of the other applications we gave for the rigidity of higher rank diagonal actions, only tells us that something is true outside an unspecified, but small, set of exceptions. The following interesting application of measure rigidity by Einsiedler, Fishman, and Shapira gives an \emph{everywhere} statement, in the spirit of Conjecture~\ref{p-adic Littlewood conjecture}:

\begin{Theorem}[Einsiedler, Fishman, and Shapira \cite{Einsiedler-Fishman-Shapira}]\label{EFS theorem} For any $\alpha \in [0, 1]$, let $n _ 1 (\alpha), n_2 (\alpha) \dots,$ be the digits in the continued fraction expansion of $\alpha$ as in \eqref{continued fractions}. Denote by $c (\alpha) = \limsup _ {i \to \infty} n _ i (\alpha)$. Then for \emph{every} irrational $\alpha \in [0, 1]$,
\begin{equation*}
\sup_ {n} \operatorname{c} \bigl(n \alpha \bmod 1\bigr) = \infty
.\end{equation*}
\end{Theorem}

Somewhat unusually, the proof of this theorem actually involves adelic dynamics \cite{Lindenstrauss-Boston}, a result closely related to Theorem~\ref{one parameter plus recurrence rigidity theorem} but with no explicit entropy assumption (the necessary entropy assumption is derived in \cite{Lindenstrauss-Boston} from the dynamical assumptions by a variation on the argument of \cite{Bourgain-Lindenstrauss}).

David Simmons observed that Theorem~\ref{EFS theorem} implies in particular that for any $\psi: \N \to \R$ with $\psi (t) \to \infty$ as $t \to \infty$, for any $\alpha \in [0, 1]$
\begin{equation} \label{towards Bourgain question}
\liminf_ {Q \to \infty} Q \min_ {q \leq Q,m \leq \psi (q)} \norm {qm \alpha} = 0
,\end{equation}
answering a question of Bourgain related to the work of Blomer, Bourgain, Radziwill and Rudnick \cite{Blomer-Bourgain-Radziwill-Rudnick} where they show that if $\alpha$ is a quadratic irrational (with some additional restrictions, removed later by Dan Carmon), for every $\epsilon > 0$, one has that $\liminf_ {Q \to \infty} Q^{2-\epsilon} \min_ {q ,m \leq Q} \norm {qm \alpha} = 0$; they also show this for a.e. $\alpha$ (but their techniques do not show \eqref{towards Bourgain question} for every $\alpha$, even when $\psi (q) = q$).

Write \[A(Q,Q') = Q \min_ {q \leq Q,m \leq Q'} \norm {qm \alpha}.\] By an (easy) result of Dirichlet $Q \min_ {q \leq Q} \norm {q \alpha}<1$ for all $\alpha, Q$, hence for any $Q,Q'$ we have the trivial estimate $A(Q,Q')\leq 1$. By considering a $\alpha \in [0,1]$ that has a sequence of extremely good approximations $\frac {p _ i }{ q _ i}$ with $q _ i$ prime, it is easy to see that there are uncountably many $\alpha$ for which $\limsup_ {Q \to \infty} A(Q,Q) =1$. However one can still gives the following strengthening of \eqref{towards Bourgain question}, whose details will appear in the forthcoming \cite{Einsiedler-Lindenstrauss-general-SL2}:

\begin{Theorem} [Einsiedler and L.] For any $\psi\to\infty$, and any $\alpha \in \R$, one has that $A(2^k,\psi(2^k)) \to 0$ outside possibly a subsequence of density zero.
\end{Theorem}

This theorem also relies on a measure classification result for higher rank diagonal actions, in this case Theorem~\ref{sl2-thm-final}.


\def\cprime{$'$}

\end{document}